\documentclass[11pt,twoside,final]{scrartcl}
\usepackage{url,a4wide}
\usepackage{amsmath, amsfonts, amssymb, amsthm, mathtools, nicefrac}
\usepackage{todonotes}
\usepackage{enumitem}
\usepackage{graphicx}
\graphicspath{{pics/}}
\usepackage{xcolor}
\usepackage[skip=1pt, font=normalsize]{subcaption}
\usepackage{hyperref,xurl} % provides \url command for bibtex and links to jump within documents
\hypersetup{plainpages=false, colorlinks, linkcolor=black, citecolor=black, urlcolor=blue,
pdftitle={On numerical realizations of Shannon's sampling theorem},
pdfauthor={Melanie Kircheis, Daniel Potts, Manfred Tasche},
pdfstartview={FitBH}}
\usepackage{placeins}
\allowdisplaybreaks

\newcommand{\e}{\mathrm e}
\renewcommand{\i}{\mathrm i}
\renewcommand{\b}{\boldsymbol}
\newcommand{\sinc}{\mathrm{sinc}}

\newcommand{\R}{\mathbb R}
\newcommand{\C}{\mathbb C}
\newcommand{\Z}{\mathbb Z}
\newcommand{\N}{\mathbb N}

\newcommand{\new}[1]{\textcolor{black}{ #1}}

\newtheorem{theorem}{Theorem}[section]
\newtheorem{corollary}[theorem]{Corollary}
\newtheorem{lemma}[theorem]{Lemma}
\newtheorem{alg}[theorem]{Algorithm}
\newtheorem{definition}[theorem]{Definition}
\newtheorem{example}[theorem]{Example}
\newtheorem{remark}[theorem]{Remark}

\newcommand{\bend}{\hspace*{0ex} \hfill \hbox{\vrule height
	1.5ex\vbox{\hrule width 1.4ex \vskip 1.4ex\hrule  width 1.4ex}\vrule
	height 1.5ex}}
\renewcommand{\qedsymbol}{\rule{1.5ex}{1.5ex}}

\newenvironment{Lemma}{\goodbreak\begin{lemma}}{\end{lemma}}
\newenvironment{Theorem}{\goodbreak\begin{theorem}}{\end{theorem}}
\newenvironment{Remark}{\goodbreak\begin{remark}\upshape}{\bend\end{remark}}

\newenvironment{Example}{\goodbreak\begin{example}\upshape}{\bend\end{example}}
\newenvironment{algorithm}[1]{\goodbreak~\begin{alg}[#1]~\vspace{-9pt}~\\
		\rule{\linewidth}{0.5pt}~\\}{\vspace{-9pt}~\\
		\rule{\linewidth}{0.5pt}~\end{alg}}

\numberwithin{equation}{section}
\numberwithin{table}{section}
\numberwithin{figure}{section}

\renewcommand{\mathbf}[1]{\ensuremath{\boldsymbol{#1}}}

\title{On numerical realizations of Shannon's sampling theorem}
\author{
	Melanie Kircheis\footnotemark[1] \and
	Daniel Potts\footnotemark[4] \and
	Manfred Tasche\footnotemark[3]
}

\date{}
\begin{document}
\maketitle

\begin{abstract}
	In this paper, we discuss some numerical realizations of Shannon's sampling theorem.
	First we show the poor convergence of classical Shannon sampling sums by presenting sharp upper and lower bounds of the norm of the Shannon sampling operator.
	In addition, it is known that in the presence of noise in the samples of a bandlimited function, the convergence of Shannon sampling series may even break down completely.
	To overcome these drawbacks, one can use oversampling and regularization with a convenient window function.
	Such a window function can be chosen either in frequency domain or in time domain.
	We especially put emphasis on the comparison of these two approaches in terms of error decay rates.
	It turns out that the best numerical results are obtained by oversampling and regularization in time domain using a $\sinh$-type window function or a continuous Kaiser--Bessel window function, which results in an interpolating
	approximation with localized sampling.
	Several numerical experiments illustrate the theoretical results.
	\medskip
	
	\smallskip
	
	AMS \emph{Subject Classifications}:
	94A20, 65T50.
\end{abstract}

\footnotetext[1]{Corresponding author: melanie.kircheis@math.tu-chemnitz.de, Chemnitz University of Technology, Faculty of Mathematics, D--09107 Chemnitz, Germany}
\footnotetext[4]{potts@mathematik.tu-chemnitz.de, Chemnitz University of
	Technology, Faculty of Mathematics, D--09107 Chemnitz, Germany}
\footnotetext[3]{manfred.tasche@uni-rostock.de, University of Rostock, Institute of Mathematics, D--18051 Rostock, Germany}

\section{Introduction}

The classical Whittaker--Kotelnikov--Shannon sampling theorem \new{(see \cite{Whittaker, Kotelnikov, Shannon49})} plays a fundamental role in signal processing, since it describes the close relation between a bandlimited function and its equidistant samples.
A function \mbox{$f \in L^2(\mathbb R)$} is called bandlimited with bandwidth~\mbox{$\frac{N}{2}$}, if
the support of its Fourier transform
\begin{equation*}
	{\hat f}(v) \coloneqq \int_{\mathbb R} f(t)\,{\mathrm e}^{-2 \pi {\mathrm i}t v}\,{\mathrm d}t\,, \quad v \in \mathbb R\,,
\end{equation*}
is contained in \mbox{$\big[- \tfrac{N}{2},\,\tfrac{N}{2}\big]$}.
Let the space of all bandlimited functions with bandwidth~\mbox{$\frac{N}{2}$} be denoted by
\begin{equation}
\label{eq:B_N/2}
	{\mathcal B}_{N/2}(\mathbb R) \coloneqq \Big\{ f \in L^2(\mathbb R)
	\colon\, \mathrm{supp}\, \hat f
	\subseteq \big[- \tfrac{N}{2},\,\tfrac{N}{2}\big] \Big\} \,,
\end{equation}
\new{which is also known as the \emph{Paley--Wiener space}.
By definition~\eqref{eq:B_N/2}, the Paley--Wiener space~\mbox{${\mathcal B}_{N/2}(\R)$} consists of equivalence classes of almost equal functions.
However, it can be shown (see e.g.~\cite[p.~6]{KiPoTa22}) that there is always a smooth representation, since for each~\mbox{$r\in\N_0$} we have by inverse Fourier transform that
\begin{align*}
	f^{(r)}(t) = \int_{-N/2}^{N/2} \hat f(v) \,(2\pi\i v)^r \,\e^{2\pi\i vt} \,\mathrm dt
\end{align*}
and \mbox{$(2\pi\i\,\cdot)^r \hat f \in L^1\big(\big[- \tfrac{N}{2},\,\tfrac{N}{2}\big]\big)$} and therefore \mbox{$f^{(r)}\in C_0(\R)$},
where \mbox{$C_0(\mathbb R)$} denotes the Banach space of continuous functions \mbox{$f \colon \R \to \C$} vanishing as \mbox{$|t| \to \infty$} with the norm
\begin{equation*}
	\|f\|_{C_0(\mathbb R)} \coloneqq \max_{t\in \mathbb R} |f(t)| \,.
\end{equation*}
That is to say, we have \mbox{${\mathcal B}_{N/2}(\R) \subseteq L^2(\R) \cap C_0(\R) \cap C^\infty(\R)$}.}

Then the Whittaker--Kotelnikov--Shannon sampling theorem states that any bandlimited function \mbox{$f \in {\mathcal B}_{N/2}(\mathbb R)$} can be recovered from its samples \mbox{$f\big(\tfrac{k}{L}\big)$},
\mbox{$k \in \mathbb Z$}, with \mbox{$L \ge N$} as
\begin{equation}
\label{eq:WKSseries}
	f(t) = \sum_{k \in \mathbb Z} f\big(\tfrac{k}{L}\big)\,\mathrm{sinc}(L\pi t - \pi k)\,, \quad t \in \mathbb R\,,
\end{equation}
where the $\sinc$ function is given by
\begin{equation*}
	\mathrm{sinc}\,x \coloneqq \left\{ \begin{array}{ll}  \frac{\sin x}{x} & \quad \colon x \in \mathbb R \setminus \{0\}\,, \\ [1ex]
		1 & \quad \colon x = 0\,. \end{array} \right.
\end{equation*}
It is well known that the series in~\eqref{eq:WKSseries} converges absolutely and uniformly on whole~$\mathbb R$.
Unfortunately, the practical use of this sampling theorem is limited, since it requires infinitely many samples, which is impossible in practice.
Furthermore, the $\sinc$ function decays very slowly such that the \emph{Shannon sampling series}
\begin{equation}
	\label{eq:Shannonseries}
	\sum_{k \in \mathbb Z} f\big(\tfrac{k}{L}\big)\,\mathrm{sinc}(L\pi t - k \pi)\,, \quad t \in \mathbb R\,,
\end{equation}
with \mbox{$L \ge N$} has rather poor convergence, as can be seen from the sharp upper and lower bounds of the norm of the Shannon sampling operator (see Theorem~\ref{Lemma:estnormST}).
In addition, it is known (see \cite{DDeV03}) that in the presence of noise in the samples \mbox{$f\big(\frac{k}{L}\big)$}, \mbox{$k \in \mathbb Z$}, of a bandlimited function \mbox{$f\in {\mathcal B}_{N/2}(\mathbb R)$}, the convergence of Shannon sampling series may even break down completely.
%In Theorem~\ref{Thm:errordata}, we show that the Shannon sampling series is not numerically robust.
%
To overcome these drawbacks, many applications employ \emph{oversampling}, i.e., a function \mbox{$f \in L^2(\mathbb R)$} of bandwidth~\mbox{$\frac{N}{2}$} is sampled on a finer grid \mbox{$\frac{1}{L}\, \mathbb Z$} with \mbox{$L > N$}, where the oversampling is measured by the \emph{oversampling parameter} \mbox{$\lambda \coloneqq \frac{L-N}{N} \ge 0$}.
In addition, we consider several regularization techniques, where a so-called \emph{window function} is used.
Since this window function can be chosen in frequency domain or in spatial domain, we study both approaches and
compare the theoretical and numerical approximation properties in terms of decay rates.

On the one hand, we investigate the regularization with a window function in frequency domain (called \emph{frequency window function}), cf. e.g. \cite{D92, Nat86, Rap96, Par97, StTa05}.
Here we use a suitable function of the form
\begin{equation*}
	{\hat \psi}(v) \coloneqq \left\{ \begin{array}{ll} 1 & \colon |v| \le \tfrac{N}{2}\,, \\
		\chi(|v|) & \colon \tfrac{N}{2}< |v| < \tfrac{L}{2}\,, \\
		0 & \colon |v| \ge \tfrac{L}{2}\,,
	\end{array} \right.
\end{equation*}
where \mbox{$\chi \colon \big[\tfrac{N}{2},\,\tfrac{L}{2}\big] \to [0,\,1]$} is frequently chosen as some monotonously decreasing, continuous function with \mbox{$\chi\big(\frac{N}{2}\big) = 1$} and \mbox{$\chi\big(\frac{L}{2}\big) = 0$}.
Applying inverse Fourier transform, we determine the corresponding function~$\psi$ in time domain. Since $\hat\psi$ is compactly supported, the uncertainty principle (cf.~\cite[p.~103, Lemma~2.39]{PPST23}) yields \mbox{$\mathrm{supp}\,\psi = \mathbb R$}.
Then it is known that the function~$f$ can be represented in the form
\begin{equation*}
	f(t) = \sum_{k \in \mathbb Z} f\big(\tfrac{k}{L}\big)\,\frac{1}{L}\,\psi\big(t - \tfrac{k}{L}\big)\,, \quad t\in \mathbb R \,.
\end{equation*}
Using uniform truncation, we approximate a function \mbox{$f \in \mathcal B_{N/2}(\mathbb R)$} by the \mbox{$T$-th} partial sum
\begin{equation*}
	(P_{\psi, T} f)(t) \coloneqq \sum_{k=-T}^T  f\big(\tfrac{k}{L}\big)\,\frac{1}{L}\,\psi\big(t - \tfrac{k}{L}\big)\,, \quad t\in [-1,\,1] \,.
\end{equation*}

On the other hand, we examine the regularization with a window function in time domain (called \emph{time window function}), cf. e.g. \cite{Q03, LZ16, MXZ09, KiPoTa22}.
Here a suitable window function \mbox{$\varphi \colon \mathbb R \to [0,\,1]$} with compact support \mbox{$\big[- \frac{m}{L},\,\frac{m}{L}\big]$} belongs to the set \mbox{$\Phi_{m/L}$} (as defined in Section~\ref{sec:time_window}) with some \mbox{$m \in \mathbb N \setminus \{1\}$}.
Then we recover a function \mbox{$f \in \mathcal B_{N/2}(\mathbb R)$} by the \emph{regularized Shannon sampling formula}
\begin{align*}
	(R_{\varphi,m} f)(t) \coloneqq \sum_{k \in \mathbb Z} f\big(\tfrac{k}{L}\big)\, {\mathrm{sinc}}\big(L \pi t - \pi k\big)\,\varphi\big(t - \tfrac{k}{L}\big) \,, \quad t\in\R \,,
\end{align*}
with \mbox{$L \ge N$}.
By defining the set \mbox{$\Phi_{m/L}$} of window functions~$\varphi$, only small truncation parameters~$m$ are needed to achieve high accuracy, resulting in short sums being computable very fast.
In other words, this approach uses localized sampling.
Moreover, this method is an interpolating approximation, since for all \mbox{$n,\,k \in \mathbb Z$} we have
\begin{align*}
	{\mathrm{sinc}}\big(L \pi t - \pi k\big)\,\varphi\big(t - \tfrac{k}{L}\big)\big|_{t=\frac nL} = \delta_{n,k}\,.
\end{align*}

In this paper we propose new estimates of the uniform approximation errors
\begin{align*}
	\max_{t \in [-1,\,1]}\big| f(t) - (P_{\psi, T} f)(t) \big|
	\quad\text{ and }\quad
	\max_{t \in \mathbb R}\big| f(t) - (R_{\varphi,m} f)(t) \big|\,,
\end{align*}
where we apply several window functions~\mbox{$\hat \psi$} and~\mbox{$\varphi$}.
\new{Note that for the frequency window functions $\hat\psi$ only error estimates on certain compact interval (such as for example \mbox{$[-1,\,1]$}) can be given, while results on whole $\mathbb R$ are not feasible.}
It is shown in the subsequent sections that the uniform approximation error decays algebraically with respect to~$T$, if \mbox{$\hat\psi$} is a frequency window function.
Otherwise, if \mbox{$\varphi \in \Phi_{m/L}$} is chosen as a time window function such as the $\sinh$-type or continuous Kaiser--Bessel window function, then the uniform approximation error decays exponentially with respect to~$m$.

To this end, this paper is organized as follows.
First, in Section~\ref{sec:Shannon} we show the poor convergence of classical Shannon sampling sums and improve results on the upper and lower bounds of the norm of the Shannon sampling operator.
Consequently, we study the different regularization techniques.
In Section~\ref{sec:frequency_window} we start with the regularization using a frequency window function.
After recapitulating a general result in Theorem~\ref{Thm:regulpsi}, we consider window functions of different regularity and present the corresponding algebraic decay results in Theorems~\ref{Thm:convregulpsi} and~\ref{Thm:err_est_psi3}.
Subsequently, in Section~\ref{sec:time_window} we proceed with the regularization using a time window function.
Here we also review the known general result in Theorem~\ref{Theorem:error-Rmf} and afterwards demonstrate the exponential decay of the considered $\sinh$-type and continuous Kaiser--Bessel window functions in Theorems~\ref{Thm:error-sinhregShannon} and \ref{Thm:error-cKBregShannon}.
Finally, in Section~\ref{sec:comp_err} we compare the previously considered approaches from Sections~\ref{sec:frequency_window} and~\ref{sec:time_window} to illustrate our theoretical results.

\section{Poor convergence of Shannon sampling sums \label{sec:Shannon}}

In order to show that the Shannon sampling series~\eqref{eq:Shannonseries} has rather poor convergence, we truncate the series~\eqref{eq:Shannonseries} with \mbox{$T\in \mathbb N$}, and consider the \mbox{$T$-\emph{th}} \emph{Shannon sampling sum}
\begin{equation}
\label{eq:ShannonSamplingSum}
	(S_T f)(t) \coloneqq \sum_{k=-T}^{T} f\big(\tfrac{k}{L}\big)\,\mathrm{sinc}(L\pi t - k \pi)\,, \quad t \in \mathbb R\,.
\end{equation}
Obviously, this operator can be formed for each \mbox{$f \in C_0(\mathbb R)$}.

\begin{Lemma}
\label{Lemma:normofST}
The linear operator \mbox{$S_T \colon \,C_0(\mathbb R) \to C_0(\mathbb R)$} has the norm
\begin{equation}
\label{eq:normofST}
	\| S_T \| = \max_{t \in \mathbb R} \sum_{k=-T}^T \big|\mathrm{sinc}(L\pi t - k \pi)\big|\,.
\end{equation}
\end{Lemma}

\emph{Proof}. For each \mbox{$f \in C_0(\mathbb R)$} and \mbox{$t\in \mathbb R$} we have
\begin{align*}
\big|(S_T f)(t) | &\le \sum_{k=-T}^{T} \big|f\big(\tfrac{k}{L}\big)\big|\,\big|\mathrm{sinc}(L\pi t - k \pi)\big| %\\
\le  \sum_{k=-T}^{T} \big|\mathrm{sinc}(L\pi t - k \pi)\big|\, \|f \|_{C_0(\mathbb R)} \,, %\le \max_{t\in \mathbb R}\sum_{k=-T}^{T} \big|\mathrm{sinc}(L\pi t - k \pi)\big|\, \|f \|_{C_0(\mathbb R)}
\end{align*}
such that
\begin{equation*}
	\| S_T f \|_{C_0(\mathbb R)} \le \max_{t\in \mathbb R}\sum_{k=-T}^{T} \big|\mathrm{sinc}(L\pi t - k \pi)\big|\, \|f \|_{C_0(\mathbb R)}\,.
\end{equation*}
By defining the even nonnegative function
\begin{equation}
\label{eq:sT}
	s_T(t) \coloneqq \sum_{k=-T}^{T} \big|\mathrm{sinc}(L\pi t - k \pi)\big|\,, \quad t \in \mathbb R\,,
\end{equation}
which is contained in \mbox{$C_0(\mathbb R)$}, and assuming that~$s_T$ has its maximum in \mbox{$t_0 \in \mathbb R$}, this yields
\begin{equation*}
	\| S_T \| = \sup \big\{ \| S_T f \|_{C_0(\mathbb R)} \colon \| f \|_{C_0(\mathbb R)}=1 \big\} \leq s_T(t_0) \,.
\end{equation*}

The other way around, we consider the linear spline \mbox{$g\in C_0(\mathbb R)$} with nodes in \mbox{$\frac{1}{L}\,\mathbb Z$}, where
\begin{equation*}
	g\big(\tfrac{k}{L}\big) = \left\{ \begin{array}{ll} \mathrm{sign}\big(\mathrm{sinc} (L \pi t_0 - k\pi)\big) & \colon k = -T\,,\,\ldots,\,T\,,\\[0.3ex]
		0 & \colon k \in \mathbb Z \setminus \{-T\,,\,\ldots,\,T\}\,.
	\end{array} \right.
\end{equation*}
Obviously, we have \mbox{$\|g\|_{C_0(\mathbb R)} = 1$} and
\begin{equation*}
	(S_T g)(t) = \sum_{k=-T}^T \mathrm{sign}\big(\mathrm{sinc} (L \pi t_0 - k\pi)\big)\, \mathrm{sinc} (L \pi t - k\pi)\big) \le s_T (t) \leq s_T(t_0) \,.
\end{equation*}
Then
\begin{equation*}
	(S_T g)(t_0) = \sum_{k=-T}^T \big| \mathrm{sinc} (L \pi t_0 - k\pi)\big| = \max_{t\in \mathbb R} s_T (t) = s_T (t_0)
\end{equation*}
implies
\begin{equation*}
	\| S_T \| \geq \| S_T g \|_{C_0(\mathbb R)} = \max_{t \in \mathbb R} \big|(S_T g)(t)\big| = s_T (t_0)
\end{equation*}
and hence~\eqref{eq:normofST}. \hfill\qedsymbol
\medskip

Now we show that the norm~\mbox{$\|S_T \|$} is unbounded with respect to~$T$.

\begin{Theorem}
\label{Lemma:estnormST}
The norm of the operator \mbox{$S_T \colon \,C_0(\mathbb R) \to C_0(\mathbb R)$} can be estimated by
\begin{equation}
\label{eq:estnormSTbylnT}
	\frac{2}{\pi}\,\big[\ln T + 2 \ln 2 + \gamma\big]- \frac{1}{\pi T\,(2T+1)} < \|S_T\| < \frac{2}{\pi}\,\big[\ln T + 2 \ln 2 + \gamma\big] + \frac{T+2}{\pi T\,(T+1)} \,,
\end{equation}
where we use \emph{Euler's constant}
\begin{equation*}
	\gamma \coloneqq \lim_{T \to \infty} \Bigg( \sum_{k=1}^T \frac{1}{k} - \ln T \Bigg) = 0.57721566 \dots \,.
\end{equation*}
\end{Theorem}

\emph{Proof}. As \new{suggested} in \cite[\new{p.~142, Problem~3.1.5}]{S93}, we represent~\new{\eqref{eq:sT}} in the form
\begin{equation*}
	s_T (t) = \sum_{k=1}^{T+1} a_k(t)\,, \quad t \in \mathbb R\,,
\end{equation*}
with
\begin{equation*}
	a_k(t) \coloneqq \left\{ \begin{array}{ll} \big|\mathrm{sinc} (L \pi t - k\pi)\big| + \big|\mathrm{sinc} (L \pi t + (k-1)\pi)\big| & \colon k = 1\,,\,\ldots,\,T\,,\\ [1ex]
	\big|\mathrm{sinc} (L \pi t + T\pi)\big| & \colon k = T+1\,.
	\end{array} \right.
\end{equation*}
Since~\eqref{eq:sT} is even, we estimate the maximum of \mbox{$s_T(t)$} only for \mbox{$t \ge 0$}.
For \mbox{$k= \new{1},\, \ldots,\,T,$} we have \new{\mbox{$a_k(0) = a_k\big(\frac 1L\big)=0$}} and
%
%For~\mbox{$a_1(t)$} with \mbox{$t \in \big(0,\, \tfrac{1}{L}\big)$} we have by trigonometric identities that
%%
%\begin{equation}
%\label{eq:a1}
%	a_1(t) = \frac{\sin (L \pi t)}{\pi}\,\bigg(\frac{1}{Lt} + \frac{1}{1-Lt}\bigg) = \frac{\sin(L \pi t)}{\pi L t\,(1-Lt)}\,.
%\end{equation}
%%
%By Schur's expansion of \mbox{$\sin (\pi x)$}, see \cite{BT88}, we know that
%%
%\begin{equation*}
%	\sin (L \pi t) = \sum_{n=1}^{\infty} \frac{1}{n!}\,\alpha_n\,(Lt)^n\,(1-Lt)^n\,, \quad t\in \big[0,\, \tfrac{1}{L}\big]\,,
%\end{equation*}
%%
%with positive coefficients
%%
%\begin{equation*}
%	\alpha_n = \frac{\pi}{(n-1)!}\,\int_0^{\pi/2} t^{n-1}\,(\pi - t)^{n-1}\,\sin t\, {\mathrm d}t\,, \quad n \in \mathbb N\,.
%\end{equation*}
%%
%Thus, we obtain the expansion
%%
%\begin{equation*}
%	a_1(t) = 1 + \frac{1}{\pi}\, \sum_{n=2}^{\infty} \frac{1}{n!}\,\alpha_n\,(Lt)^{n-1}\,(1-Lt)^{n-1}\,, \quad t\in \big[0,\, \tfrac{1}{L}\big]\,.
%\end{equation*}
%%
%Hence, the function \mbox{$a_1 \colon \,\big[0,\, \tfrac{1}{L}\big] \to \mathbb R$} is concave and has its maximum at \mbox{$t = \frac{1}{2L}$}, i.e., by~\eqref{eq:a1} we compute
%%
%\begin{equation*}
%	\max_{t \in [0,\,1/L]} a_1(t) = a_1\big(\tfrac{1}{2L}\big) = \frac{4}{\pi}\,.
%\end{equation*}
%%
%For \mbox{$a_k(t)$}, \mbox{$k= 2,\, \ldots,\,T$}, with \mbox{$t \in \big(0,\, \tfrac{1}{L}\big)$} we have
by trigonometric identities we obtain for \mbox{$t \in \big(0,\, \tfrac{1}{L}\big)$} that
\begin{equation*}
	a_k(t)
	= \frac{\sin (L \pi t)}{\pi}\,\bigg(\frac{1}{k-1+Lt} + \frac{1}{k-Lt}\bigg)
	= \frac{(2k-1)\,\sin(L \pi t)}{\pi\,\big[(k-1)k +Lt\,(1-Lt)\big]}\,.
\end{equation*}
We define the functions \mbox{$b_k \colon \,(0,\,1) \to \mathbb R$}, \mbox{$k = \new{1},\, \ldots,\,T$}, via
\begin{equation*}
	b_k(x) \coloneqq \new{\frac{\pi}{(2k-1)}} \, a_k\big(\tfrac{x}{L}\big) = \frac{\sin (\pi x)}{(k-1)k + x (1-x)}\,, \quad x \in (0,\, 1) \,,
\end{equation*}
such that the symmetry relation \mbox{$b_k(x) = b_k(1-x)$} is fulfilled, i.e., each~\mbox{$b_k$} is symmetric with reference to~\mbox{$\frac{1}{2}$}.
Furthermore, by \mbox{$b_k'(x) \geq 0$} for \mbox{$x \in \new{\big(}0,\,\frac{1}{2}\big]$}, the function~$b_k$ is increasing on \mbox{$\new{\big(}0,\,\frac{1}{2}\big]$} and therefore has its maximum at \mbox{$x = \tfrac{1}{2}$}.
Thus, the function \mbox{$a_k \colon \,\big[0,\, \tfrac{1}{L}\big]\to \mathbb R$} has its
maximum at \mbox{$t = \tfrac{1}{2L}$}, i.e.,
\begin{equation*}
	\max_{t \in [0,\,1/L]} a_k(t) = a_k\big(\tfrac{1}{2L}\big) = \frac{4}{(2k-1)\pi}\,.
\end{equation*}
Since \mbox{$a_{T+1}(t)$} can be written as
\begin{equation*}
	a_{T+1}(t) = \frac{\sin (L\pi t)}{\pi \,(T + Lt)}
\end{equation*}
for \mbox{$t\in \big[0,\, \tfrac{1}{L}\big]$}, we obtain
\begin{equation*}
	0 < \max_{t \in [0,\,1/L]} a_{T+1}(t) < \frac{1}{\pi T}\,.
\end{equation*}
%
%In the case \mbox{$T \gg 1$}, the function \mbox{$a_{T+1} \colon \,\big[0,\, \tfrac{1}{L}\big]\to \mathbb R$} has its maximum close to \mbox{$t = \tfrac{1}{2L}$}.
Hence, in summary this yields
\begin{equation*}
	\frac{4}{\pi}\,\sum_{k=1}^T \frac{1}{2k-1} < \max_{t \in [0,\,1/L]} s_T(t) < \frac{4}{\pi}\,\sum_{k=1}^T \frac{1}{2k-1} + \frac{1}{\pi T}\,.
\end{equation*}
For \mbox{$t\in \big[\tfrac{n}{L},\,\tfrac{n+1}{L}\big]$} with arbitrary \mbox{$n \in \mathbb N$}, the sum~\mbox{$s_T(t)$} is less than it is for \mbox{$t\in \big[0,\, \tfrac{1}{L}\big]$}, since
for each \mbox{$n\in \mathbb N$} and \mbox{$t\in \big(0, \, \tfrac{1}{L}\big)$} we have
\begin{equation*}
	\sum_{k=-T}^T \frac{\sin(L\pi t)}{|L \pi t - (k-n)\pi|} < \sum_{k=-T}^T \frac{\sin(L\pi t)}{|L \pi t - k\pi|}
\end{equation*}
and therefore \mbox{$s_T\big(\tfrac{n}{L} + t \big) < s_T(t)$}.
Thus, we obtain
%
%\begin{equation}
%\label{eq:estsT}
%	\frac{4}{\pi}\,\sum_{k=1}^T \frac{1}{2k-1} < \max_{t \in \mathbb R} s_T(t) < \frac{4}{\pi}\,\sum_{k=1}^T \frac{1}{2k-1} + \frac{1}{\pi T}\,.
%\end{equation}
%%
%By Lemma~\ref{Lemma:normofST} this can also be written as
%%
\begin{equation}
\label{eq:est1normST}
	\frac{4}{\pi}\,\sum_{k=1}^T \frac{1}{2k-1} < \max_{t \in \mathbb R} s_T(t) = \| S_T \| < \frac{4}{\pi}\,\sum_{k=1}^T \frac{1}{2k-1} + \frac{1}{\pi T}\,.
\end{equation}
Note that for \mbox{$T \gg 1$} the value
\begin{equation}
\label{eq:sT(1/(2L))}
	s_T\big(\tfrac{1}{2L}\big) = \frac{4}{\pi}\,\sum_{k=1}^T \frac{1}{2k-1} + \frac{2}{\pi\,(2 T + 1)}
\end{equation}
is a good approximation of the norm~\mbox{$\|S_T\|$}.
\medskip

Now we estimate~\mbox{$\|S_T\|$} by \mbox{$\ln T$}.
For this purpose we denote the \mbox{$T$-\emph{th}} \emph{harmonic number} by
\begin{equation*}
	H_T \coloneqq \sum_{k=1}^T \frac{1}{k}\,, \quad T \in \mathbb N\,,
\end{equation*}
such that
\begin{equation}
\label{eq:H2T-1/2HT}
	\sum_{k=1}^T \frac{1}{2k-1}
	= \sum_{k=1}^T \bigg(\frac{1}{2k-1} + \frac{1}{2k}\bigg) - \sum_{k=1}^T \frac{1}{2k}
	= H_{2T} - \frac{1}{2}\,H_T\,.
\end{equation}
Using \emph{Euler's constant}
\begin{equation*}
	\gamma = \lim_{T \to \infty} \big( H_T - \ln T\big) \,,
\end{equation*}
the estimates
\begin{equation}
\label{eq:estHT}
	\frac{1}{2T + 2}< H_T - \ln T - \gamma < \frac{1}{2T}
\end{equation}
are known (see \cite{Y91}).
From~\eqref{eq:H2T-1/2HT} and~\eqref{eq:estHT} we conclude that
\begin{equation}
\label{eq:estH2T-1/2HT}
	\frac{1}{2}\,\ln T + \ln 2 + \frac{1}{2}\,\gamma - \frac{1}{4T\,(2T+1)} < H_{2T}- \frac{1}{2}\,H_T < \frac{1}{2}\,\ln T + \ln 2 + \frac{1}{2}\,\gamma + \frac{1}{4 T\,(T+1)}\,.
\end{equation}
Therefore, applying~\eqref{eq:est1normST}, \eqref{eq:H2T-1/2HT}, and~\eqref{eq:estH2T-1/2HT} yields the assertion~\eqref{eq:estnormSTbylnT}. \hfill\qedsymbol
\medskip

\begin{figure}[ht]
	\centering
	\captionsetup[subfigure]{justification=centering}
	\begin{subfigure}[t]{0.435\textwidth}
		\includegraphics[width=\textwidth]{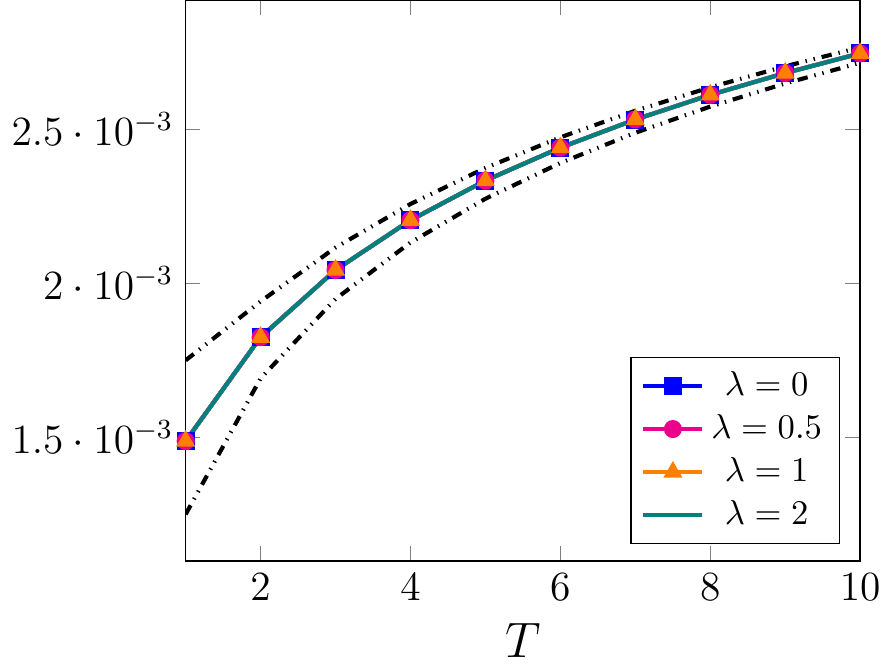}
%		\caption{old version}
	\end{subfigure}
	\begin{subfigure}[t]{0.4\textwidth}
		\includegraphics[width=\textwidth]{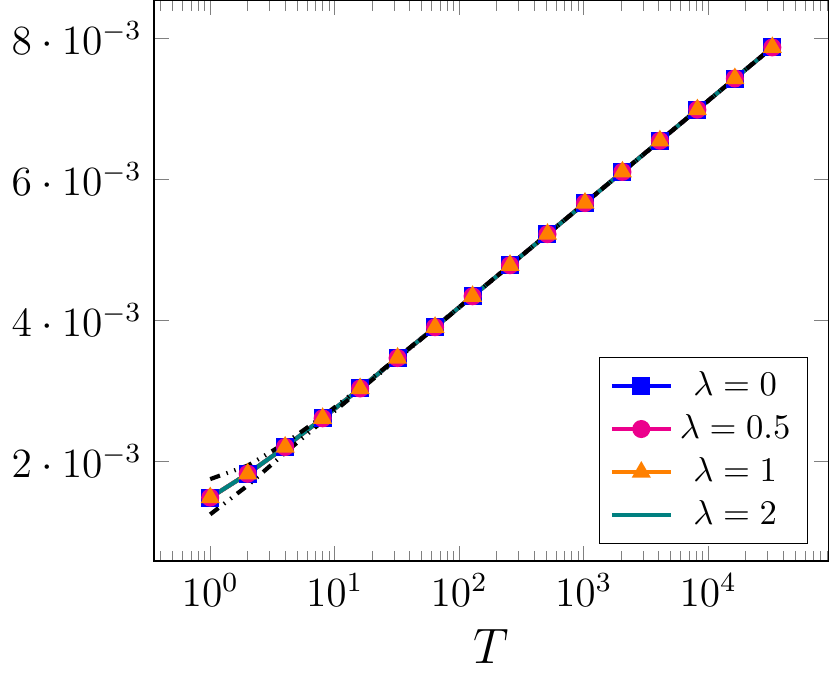}
%		\caption{new version}
	\end{subfigure}
	\caption{The norm \mbox{$\|\tilde f - f \|_{C_0(\mathbb R)}$} of~\eqref{eq:ftilde-f} as well as its lower/upper bounds~\eqref{eq:lower_bound} and~\eqref{eq:upper_bound} for several~\new{\mbox{$T \in \N$}} and~\mbox{$L=N(1+\lambda)$} with \mbox{$\lambda \in \{0,0.5,1,2\}$}, where~\mbox{$N=128$} and~\mbox{$\varepsilon=10^{-3}$} are chosen.
		\label{fig:nonrobustness}}
\end{figure}

\new{Note that Theorem~\ref{Lemma:estnormST} improves a former result of \cite[p.~142, Problem 3.1.5]{S93}, which only contains a coarse proof sketch for the upper bound, while~\eqref{eq:estnormSTbylnT} gives a very precise nesting, see also Figure~\ref{fig:nonrobustness}.}
Additionally, we remark that \new{Theorem}~\ref{Lemma:estnormST} immediately implies
\begin{equation*}
	\lim_{T \to \infty} \bigg( \|S_T\| - \frac{2}{\pi}\,\ln T \bigg) = \frac{4}{\pi}\,\ln 2 + \frac{2 \gamma}{\pi}\,.
\end{equation*}
%
%\medskip

\begin{Remark}
\new{
By Theorem~\ref{Lemma:estnormST} and the theorem of Banach--Steinhaus, an arbitrary function~\mbox{$f \in C_0(\R)$} cannot be represented in the form~\eqref{eq:Shannonseries}, since the norm of the linear operator \mbox{$S_T \colon \,C_0(\R) \to C_0(\R)$} is unbounded with respect to~$T$.
However, as known from the Whittaker--Kotelnikov--Shannon sampling theorem for bandlimited functions~\mbox{$f \in {\mathcal B}_{N/2}(\R)$},
the series~\eqref{eq:Shannonseries} converges absolutely and uniformly on whole~\mbox{$\R$}. \\
\indent
Nevertheless, this convergence is very slow due to the poor decay of the \mbox{$\sinc$ func}tion, as can be seen from the sharp upper and lower bounds of the norm of the Shannon sampling operator, see Lemma~\ref{Lemma:estnormST}.
More precisely, rigorous analysis of the approximation error, cf.~\cite[Theorem~1]{Ja66}, shows that on the fixed interval \mbox{$[-1,\,1]$} we have a convergence rate of only~\mbox{$(T-L)^{-1/2}$},
which was also mentioned in~\cite{Za93, StTa06}. \\
\indent
Note that general results on whole~$\R$ are not feasible for arbitrary~\mbox{$f \in {\mathcal B}_{N/2}(\R)$}.
Rather, the uniform approximation error~\mbox{$\| f - S_T f\|_{C_0(\R)}$} can only be studied under the additional assumption that~\mbox{$f \in {\mathcal B}_{N/2}(\R)$} satisfies certain decay conditions, see~\cite{Li98, JiGe03}.}
\end{Remark}

Now let \mbox{$f \in {\mathcal B}_{N/2}(\mathbb R)$} be a given bandlimited function with bandwidth \mbox{$\tfrac{N}{2}$} and let \mbox{$T \in \mathbb N$} be sufficiently large.
For given samples \mbox{$f\big(\tfrac{k}{L}\big)$} with \mbox{$k\in \mathbb Z$} and \mbox{$L\ge N$} we consider finitely many erroneous samples
\begin{equation*}
	{\tilde f}_k \coloneqq  \left\{ \begin{array}{ll} f\big(\tfrac{k}{L}\big) + \varepsilon_k & \colon k = -T,\,\ldots,\,T\,,\\ [1ex]
	f\big(\tfrac{k}{L}\big) & \colon k \in \mathbb Z \setminus \{-T,\,\ldots,\,T\} \,,
	\end{array} \right.
\end{equation*}
with error terms \mbox{$\varepsilon_k$} which are bounded by \mbox{$|\varepsilon_k| \le \varepsilon$} for \mbox{$k = -T,\, \ldots,\,T$}.
Then for the approximation
\begin{equation*}
	{\tilde f}(t) \coloneqq \sum_{k\in \mathbb Z} {\tilde f}_k\, \mathrm{sinc}(L\pi t - k \pi) = f(t) + \sum_{k=-T}^T \varepsilon_k\,\mathrm{sinc}(L\pi t - k \pi)\,, \quad t\in \mathbb R \,,
\end{equation*}
\new{the following error bounds can be shown.}

\begin{Theorem}
\label{Thm:errordata}
Let \mbox{$f \in {\mathcal B}_{N/2}(\mathbb R)$} be an arbitrary bandlimited function with bandwidth \mbox{$\tfrac{N}{2}$}.
Further let \mbox{$L \ge N$}, \mbox{$T \in \mathbb N$}, and \mbox{$\varepsilon > 0$} be given.\\
\new{Then we have}
\begin{align}
\label{eq:upper_bound}
	\new{\|\tilde f - f\|_{C_0(\mathbb R)}
	< \varepsilon\,\bigg(\frac{2}{\pi}\,\ln T + \frac{5}{4} + \frac{1}{2T}\bigg) \,.}
\end{align}
Moreover, for the special error terms
\begin{align*}
	\varepsilon_k &= \varepsilon\, \mathrm{sign}\big(\mathrm{sinc}(\tfrac{\pi}{2} - k \pi)\big) %\\
	= \varepsilon\,(-1)^{k+1}\,{\mathrm{sign}}(2k-1)\,, \quad k = -T,\ldots,T\,,
\end{align*}
we have
\begin{equation}
\label{eq:lower_bound}
	\|{\tilde f} - f \|_{C_0(\mathbb R)} \ge \varepsilon\,\bigg(\frac{2}{\pi}\,\ln T + \frac{4}{\pi}\,\ln 2 + \frac{2 \gamma}{\pi}\bigg) > \varepsilon\,\bigg(\frac{2}{\pi}\,\ln T + \frac{5}{4}\bigg)\,,
\end{equation}
such that the Shannon sampling series is not numerically robust in the worst case analysis.
\end{Theorem}

\emph{Proof}.
\new{By~\eqref{eq:sT}, \eqref{eq:est1normST} and~\eqref{eq:estH2T-1/2HT} we are given the upper error bound
\begin{align*}
	\|\tilde f - f\|_{C_0(\mathbb R)}
	\le \varepsilon\,\max_{t\in \mathbb R} s_T(t)
	< \frac{4 \varepsilon}{\pi}\,\sum_{k=1}^T \frac{1}{2k-1} + \frac{\varepsilon}{\pi T}
	< \frac{2\varepsilon}{\pi}\, \bigg(\ln T + 2\ln 2 + \gamma + \frac{T+2}{2T(T+1)}\bigg) \,.
\end{align*}
Since
\begin{equation}
\label{eq:gamma_term}
	\frac{4}{\pi}\,\ln 2 + \frac{2 \gamma}{\pi} = 1.2500093\ldots
\end{equation}
and \mbox{$\mu_T \coloneqq \frac{T+2}{T+1} = 1+\frac{1}{T+1}$} is monotonously decreasing with \mbox{$\max_{T\in \mathbb N} \mu_T = \mu_1 = \tfrac 32$}, we have
\begin{align*}
	\frac{4}{\pi}\,\ln 2 + \frac{2 \gamma}{\pi} + \frac{1}{\pi T}\cdot\mu_T
	\leq
	\frac 54 + \frac{1}{2T} \,,
\end{align*}
which yields~\eqref{eq:upper_bound}.}

Next, we proceed with the lower bound.
Due to the special choice of the error terms~\mbox{$\varepsilon_k$} we obtain
\begin{equation}
\label{eq:ftilde-f}
	{\tilde f}(t) - f(t) = \varepsilon\, \sum_{k=-T}^T \mathrm{sign}\big(\mathrm{sinc}(\tfrac{\pi}{2} - k \pi)\big)\;\mathrm{sinc}(L \pi t - k \pi)\,, \quad t\in \mathbb R\,.
\end{equation}
By~\eqref{eq:sT(1/(2L))} and~\eqref{eq:estH2T-1/2HT} we conclude that
\begin{align*}
	\|{\tilde f} - f \|_{C_0(\mathbb R)} &\ge \big|{\tilde f}(\tfrac{1}{2L}) - f(\tfrac{1}{2L})\big| = \varepsilon\, \sum_{k=-T}^T \big|\mathrm{sinc}(\tfrac{\pi}{2} - k \pi)\big|
	= \varepsilon\,s_T\big(\tfrac{1}{2L}\big) \\[1ex]
	&= \frac{4\varepsilon}{\pi}\,\sum_{k=1}^T \frac{1}{2k-1} + \frac{2 \varepsilon}{(2T+1)\pi} %\\[1ex]
%	> \varepsilon\,\bigg(\frac{2}{\pi}\,\ln T + \frac{4}{\pi}\,\ln 2 + \frac{2 \gamma}{\pi}\bigg) + \varepsilon\,\frac{2T-1}{(2T+1)T\,\pi} \\[1ex]
	> \varepsilon\,\bigg(\frac{2}{\pi}\,\ln T + \frac{4}{\pi}\,\ln 2 + \frac{2 \gamma}{\pi}\bigg) \,,
\end{align*}
such that~\eqref{eq:gamma_term} completes the proof. \hfill\qedsymbol
\medskip

\new{Note that~\eqref{eq:upper_bound} specifies a corresponding remark of \cite[p.~681]{DDeV03} as it makes the constant explicit.}
In addition, we remark that the norm \mbox{$\|\tilde f - f \|_{C_0(\mathbb R)}$} does not depend on the special choice of the function~$f$ or the oversampling parameter~$\lambda$,
see Figure~\ref{fig:nonrobustness}.
Furthermore, Figure~\ref{fig:nonrobustness} also illustrates that \new{for \mbox{$T \to \infty$}} the error behavior shown in Theorem~\ref{Thm:errordata} is not satisfactory.
Thus, in the presence of noise in the samples~\mbox{$f\big(\tfrac{{k}}{L}\big)$}, \mbox{${k} \in \Z$}, the convergence of the Shannon sampling series~\eqref{eq:Shannonseries} may even break down completely.

\begin{Remark}
\label{Rem:average case}
\new{In the above worst case analysis we have seen that the approximation of \mbox{$f \in {\mathcal B}_{N/2}(\mathbb R)$} by the \mbox{$T$-th} partial sum~\eqref{eq:ShannonSamplingSum} of its Shannon sampling series
%\begin{equation}
%\label{eq:T-thpartsum}
%\sum_{k=-T}^T f\big(\tfrac{k}{L}\big)\,\mathrm{sinc}(L\pi t - k \pi)
%\end{equation}
with~\mbox{$L \geq N$} is not numerically robust in the deterministic sense. Otherwise, a simple average case study (see~\cite{TZ01}) shows that this approximation is numerically robust  in the stochastic sense.
Therefore, we compute~\eqref{eq:ShannonSamplingSum} as an inner product of the real \mbox{$(2T+1)$-dim}ensional vectors
$$
\Big( f\big(\tfrac{k}{L}\big) \Big)_{k=-T}^T \quad\text{ and }\quad \big( \mathrm{sinc}(L\pi t - k \pi) \big)_{k=-T}^T\,.
$$
Now assume that instead of the exact samples~\mbox{$f\big(\frac{k}{L}\big)$}, \mbox{$k = - T, \ldots,T$}, only perturbed samples \mbox{${\tilde f}\big(\frac{k}{L}\big) =  f\big(\frac{k}{L}\big) + X_k$}, \mbox{$k = - T, \ldots,T$}, are given, where
the real random variables~$X_k$ are uncorrelated, each having expectation \mbox{${\mathbb E}(X_k) = 0$} and constant variance \mbox{${\mathbb V}(X_k) = {\mathbb E}(|X_k|^2) = \rho^2$} with~\mbox{$\rho > 0$}.
Then we consider the stochastic approximation error
\begin{eqnarray*}
\Delta_T &:=& \sum_{k=-T}^T {\tilde f}\big(\tfrac{k}{L}\big)\,\mathrm{sinc}(L\pi t - k \pi) - \sum_{k=-T}^T f\big(\tfrac{k}{L}\big)\,\mathrm{sinc}(L\pi t - k \pi)\\
&=& \sum_{k=-T}^T X_k\,\mathrm{sinc}(L\pi t - k \pi)\,.
\end{eqnarray*}
Obviously, this error term~$\Delta_T$ has the expectation
$$
{\mathbb E}(\Delta_T) = \sum_{k=-T}^T \mathrm{sinc}(L\pi t - k \pi)\,{\mathbb E}(X_k) = 0
$$
and the variance
$$
{\mathbb V}(\Delta_T) = \sum_{k=-T}^T \big|\mathrm{sinc}(L\pi t - k \pi)\big|^2\,{\mathbb V}(X_k) = \rho^2\,\sum_{k=-T}^T \big|\mathrm{sinc}(L\pi t - k \pi)\big|^2\,.
$$
From \cite[p.~89, Problem 1.10.9]{S93} it follows that
$$
\sum_{k\in \mathbb Z} \big|\mathrm{sinc}(L\pi t - k \pi)\big|^2 = 1 \,,
$$
such that \mbox{${\mathbb V}(\Delta_T) \leq \rho^2$}.}
\end{Remark}

\section{Regularization with a frequency window function \label{sec:frequency_window}}

To overcome the drawbacks of poor convergence and numerical instability, one can apply regularization with a convenient window function either in the frequency domain or in the time domain.
Often one employs \emph{oversampling}, i.e., a bandlimited function \mbox{$f \in {\mathcal B}_{N/2}(\mathbb R)$} of bandwidth~$\tfrac N2$ is sampled on a finer grid \mbox{$\frac{1}{L}\, \mathbb Z$} with \mbox{$L =N\,(1+\lambda)$}, where \mbox{$\lambda >0$} is the oversampling parameter.
\medskip

First, together with oversampling, we consider the  \emph{regularization with a frequency window function} of the form
\begin{equation}
\label{eq:hatpsi(v)}
	{\hat \psi}(v) \coloneqq \left\{ \begin{array}{ll} 1 & \colon |v| \le \tfrac{N}{2}\,, \\
	\chi(|v|) & \colon \tfrac{N}{2}< |v| < \tfrac{L}{2}\,, \\
	0 & \colon |v| \ge \tfrac{L}{2}\,,
	\end{array} \right.
\end{equation}
cf.~\cite{D92, Par97, StTa05}, where \mbox{$\chi \colon \big[\tfrac{N}{2},\,\tfrac{L}{2}\big] \to [0,\,1]$} is frequently chosen as some monotonously decreasing, continuous function with \mbox{$\chi\big(\frac{N}{2}\big) = 1$} and \mbox{$\chi\big(\frac{L}{2}\big) = 0$}.
Applying the inverse Fourier transform, we determine the corresponding function in time domain as
\begin{equation}
\label{eq:psi}
	\psi(t) = \int_{\mathbb R} {\hat \psi}(v)\, {\mathrm e}^{2 \pi {\mathrm i}\,v t}\, {\mathrm d}v = 2\,\int_0^{L/2} {\hat \psi}(v)\, \cos(2 \pi\,v t)\, {\mathrm d}v\,.
\end{equation}
\begin{Example}
\label{Ex:hatpsilin}
A simple example of a frequency window function is the  \emph{linear frequency window function} (cf. \cite[pp.~18--19]{D92} or \cite[pp.~210--212]{Par97})
\begin{equation}
\label{eq:hatpsilin(v)}
	{\hat \psi_{\mathrm{lin}}}(v)
	\coloneqq \left\{ \begin{array}{ll} 1 & \colon |v| \le \tfrac{N}{2}\,, \\
		1 - \tfrac{2 |v| - N}{L-N} & \colon \tfrac{N}{2}< |v| < \tfrac{L}{2}\,, \\
		0 & \colon |v| \ge \tfrac{L}{2}\,.
	\end{array} \right.
\end{equation}
\new{Note that in a trigonometric setting a function of the form~\eqref{eq:hatpsilin(v)} is also often referred to as trapezoidal or de La Vall\'{e}e-Poussin type window function, respectively.}
Obviously, \mbox{${\hat \psi_{\mathrm{lin}}}(v)$} is a continuous linear spline supported on \mbox{$\big[-\tfrac{L}{2},\,\tfrac{L}{2}\big]$}, see Figure~\ref{fig:frequency_window}~(a).
By~\eqref{eq:psi} we receive \mbox{$\psi_{\mathrm{lin}}(0) = \tfrac{N+L}{2}$}.
For \mbox{$t\in \mathbb R \setminus \{0\}$} we obtain
\new{\begin{align}
\label{eq:psilin}
	\psi_{\mathrm{lin}}(t) &= 2\,\int_0^{N/2} \cos(2 \pi\,v t)\, {\mathrm d}v + 2\,\int_{N/2}^{L/2} \bigg(1 - \frac{2v-N}{L-N}\bigg)\, \cos(2 \pi\,v t)\, {\mathrm d}v \notag \\[1ex]
%	&= \frac{1}{(L-N)\,(\pi t)^2}\,\big(\cos(N\pi t) - \cos(L\pi t)\big) \notag \\[1ex]
%	&= \frac{2}{(L-N)\,(\pi t)^2}\, \sin\big(\tfrac{N+L}{2}\,\pi t\big) \,\sin\big(\tfrac{L-N}{2}\,\pi t\big) \notag \\[1ex]
	&= \frac{N+L}{2}\,\mathrm{sinc}\big(\tfrac{N+L}{2}\,\pi t\big) \,\mathrm{sinc}\big(\tfrac{L-N}{2}\,\pi t\big)\,.
\end{align}}
This function \mbox{$\tfrac{1}{L}\, \psi_{\mathrm{lin}}$} is even, supported on whole~$\mathbb R$, has its maximum at \mbox{$t=0$} such that
\begin{equation*}
	\Big\|\frac{1}{L}\, \psi_{\mathrm{lin}} \Big\|_{C_0(\mathbb R)} = \frac{1}{L}\, \psi_{\mathrm{lin}}(0) = \frac{2 + \lambda}{2 + 2\lambda} < 1 \,.
\end{equation*}
In addition, \mbox{$\tfrac{1}{L}\, \psi_{\mathrm{lin}}(t)$} has a faster decay than \mbox{$\mathrm{sinc}(N \pi t)$} for \mbox{$|t| \to \infty$}, cf.~Figure~\ref{fig:frequency_window}~(b).
Note that we have
\begin{equation*}
	\lim_{L \to +N} \frac{1}{L}\, \psi_{\mathrm{lin}}(t) = \mathrm{sinc}(N \pi t)\,.
\end{equation*}
\begin{figure}[ht]
	\centering
	\captionsetup[subfigure]{justification=centering}
	\begin{subfigure}[t]{0.4\textwidth}
		\includegraphics[width=\textwidth]{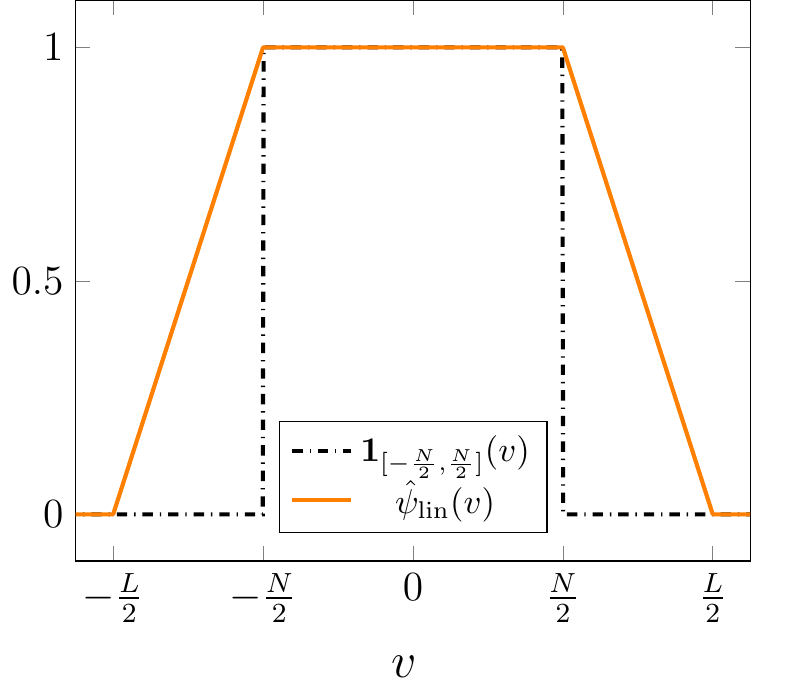}
		\caption{$\hat\psi_{\mathrm{lin}}$ in~\eqref{eq:hatpsilin(v)}}
	\end{subfigure}
	\begin{subfigure}[t]{0.4\textwidth}
		\includegraphics[width=\textwidth]{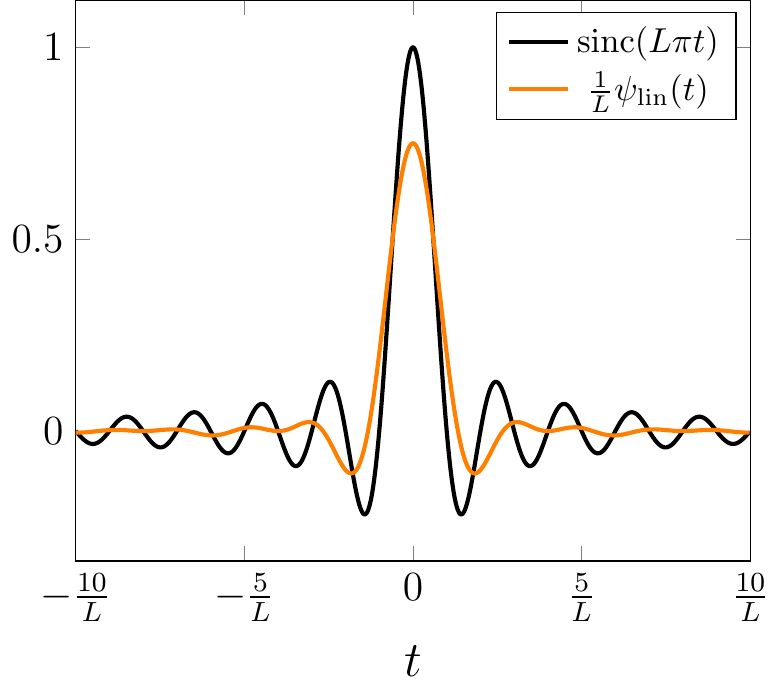}
		\caption{$\frac{1}{L}\,\psi_{\mathrm{lin}}$ in~\eqref{eq:psilin}}
	\end{subfigure}
	\caption{The frequency window function~\eqref{eq:hatpsilin(v)} and its scaled inverse Fourier transform~\eqref{eq:psilin}.
	\label{fig:frequency_window}}
\end{figure}
\end{Example}

\begin{Remark}
Note that \mbox{$\big\{\tfrac{1}{L}\,\psi_{\mathrm{lin}}\big(\cdot\, - \tfrac{k}{L}\big)\big\}_{k \in \mathbb Z}$} is a \emph{Bessel sequence in \mbox{$L^2(\mathbb R)$} with the bound} \mbox{$\tfrac{1}{L}$}, i.e., for all \mbox{$f \in L^2(\mathbb R)$} we have
\begin{equation*}
	\sum_{k \in \mathbb Z} \big| \langle f,\,\tfrac{1}{L}\,\psi_{\mathrm{lin}}\big(\cdot\, - \tfrac{k}{L}\big) \rangle_{L^2(\mathbb R)}\big|^2 \le \tfrac{1}{L}\,\|f\|_{L^2(\mathbb R)}^2\,.
\end{equation*}
However, \mbox{$\big\{\tfrac{1}{L}\,\psi_{\mathrm{lin}}\big(\cdot\, - \tfrac{k}{L}\big)\big\}_{k \in \mathbb Z}$} is not an orthonormal sequence and also not a Riesz sequence.
To see this, we consider the 1-periodic function
\begin{equation*}
	\Psi(v) \coloneqq \frac{1}{L^2}\,\sum_{k \in \mathbb Z}\big| {\hat \psi_{\mathrm{lin}}}(L v + L k)\big|^2\,, \quad v \in \mathbb R\,.
\end{equation*}
By~\eqref{eq:hatpsilin(v)} we have
\begin{equation*}
	{\hat \psi_{\mathrm{lin}}}(Lv) = \left\{ \begin{array}{ll} 1 & \colon |v| \le \tfrac{1}{2 + 2\lambda}\,, \\[1ex]
	1 - \tfrac{(2 + 2\lambda) |v| - 1}{\lambda} & \colon \tfrac{1}{2+2\lambda}< |v| < \tfrac{1}{2}\,, \\[1ex]
	0 & \colon |v| \ge \tfrac{1}{2}\,,
	\end{array} \right.
\end{equation*}
i.e., \mbox{$\Psi(v) \le \frac{1}{L^2}$} for all \mbox{$v \in \mathbb R$} and \mbox{$\Psi(v) = 0$} for \mbox{$v \in \tfrac{1}{2} + \mathbb Z$}.
Then \cite[Theorem~9.2.5]{Ch16} yields the result.
\end{Remark}
\medskip

Analogous to \cite[p.~19]{D92} and \cite[Theorem~7.2.5]{Par97}, we obtain the following representation result, see also \cite[p.~4]{StTa05}.

\begin{Theorem}
\label{Thm:regulpsi}
Let \mbox{$f \in {\mathcal B}_{N/2}(\mathbb R)$} be a bandlimited function with bandwidth \mbox{$\tfrac{N}{2}$}, \mbox{$N \in \N$}, and let \mbox{$L = N\,(1+\lambda)$} with \mbox{$\lambda > 0$} be given.
\new{Assume that the samples~\mbox{$f\big(\frac{k}{L}\big)$}, \mbox{$k \in \mathbb Z$}, fulfill the condition}
\begin{equation}
\label{eq:sampledecay}
\new{\sum_{k\in \mathbb Z} \big| f\big(\tfrac{k}{L}\big) \big| < \infty\,.}
\end{equation}
Using oversampling and regularization
with a frequency window function~\mbox{$\hat \psi$} of the form~\eqref{eq:hatpsi(v)}, the function~$f$ can be represented as
\begin{equation}
\label{eq:f=sum}
	f(t) = \sum_{k \in \mathbb Z} f\big(\tfrac{k}{L}\big)\,\frac{1}{L}\,\psi\big(t - \tfrac{k}{L}\big)\,, \quad t\in \mathbb R \,,
\end{equation}
\new{where the series~\eqref{eq:f=sum} converges absolutely and uniformly on $\mathbb R$.}
\end{Theorem}

\emph{Proof}. Since by assumption \mbox{$f \in {\mathcal B}_{N/2}(\mathbb R)$}, we have \mbox{$\mathrm{supp}\,\hat f \subseteq \big[- \tfrac{N}{2},\,\tfrac{N}{2}\big] \subset \big[- \tfrac{L}{2},\,\tfrac{L}{2}\big]$} and therefore the function~${\hat f}$ restricted on \mbox{$\big[- \tfrac{L}{2},\,\tfrac{L}{2}\big]$}
can be represented by its \mbox{$L$-periodic} Fourier series
\begin{equation}
\label{eq:Fourierserieshatf}
	{\hat f}(v) = \sum_{k \in \mathbb Z} c_k(\hat f)\, {\mathrm e}^{2 \pi {\mathrm i}\,k v/L}\,, \quad v \in \big[- \tfrac{L}{2},\,\tfrac{L}{2}\big]\,,
\end{equation}
with the Fourier coefficients
\begin{equation*}
	c_k(\hat f) = \frac{1}{L}\, \int_{-L/2}^{L/2} {\hat f}(v)\,{\mathrm e}^{-2 \pi {\mathrm i}\,k v/L}\,{\mathrm d}v\,.
\end{equation*}
Using the inverse Fourier transform, we see that
\begin{equation*}
	c_k(\hat f) = \frac{1}{L}\, \int_{\mathbb R} {\hat f}(v)\,{\mathrm e}^{-2 \pi {\mathrm i}\,k v/L}\,{\mathrm d}v = \frac{1}{L}\,f\big(-\tfrac{k}{L}\big)\,.
\end{equation*}
Hence, we may write~$\hat f$ given by~\eqref{eq:Fourierserieshatf} as
\begin{equation}
\label{eq:Fourierserieshatfnew}
	{\hat f}(v) =  \frac{1}{L}\,\sum_{k \in \mathbb Z}  f\big(\tfrac{k}{L}\big) \, {\mathrm e}^{-2 \pi {\mathrm i}\,k v/L}\,, \quad v \in \big[- \tfrac{L}{2},\,\tfrac{L}{2}\big]\,.
\end{equation}
\new{By the Weierstrass M-test and assumption~\eqref{eq:sampledecay}, the Fourier series~\eqref{eq:Fourierserieshatfnew} converges absolutely and uniformly on \mbox{$\big[- \tfrac{L}{2},\,\tfrac{L}{2}\big]$}.}
Additionally, we have \mbox{${\hat \psi}(v) = 1$} for \mbox{$v \in \big[- \tfrac{N}{2},\, \tfrac{N}{2}\big]$} by~\eqref{eq:hatpsi(v)} as well as \mbox{$\mathrm{supp}\,\hat f \subseteq \big[- \tfrac{N}{2},\,\tfrac{N}{2}\big]$} by assumption, such that \mbox{${\hat f}(v) = {\hat f}(v)\,{\hat \psi}(v)$} for all \mbox{$v \in \mathbb R$}. Therefore, we obtain
\begin{align*}
	f(t)
	&= \int_{\mathbb R} {\hat f}(v)\,{\mathrm e}^{2 \pi {\mathrm i}\,t v}\,{\mathrm d}v = \int_{-L/2}^{L/2} {\hat f}(v)\,{\mathrm e}^{2 \pi {\mathrm i}\,t v}\,{\mathrm d}v %\\[1ex]
	= \int_{-L/2}^{L/2} {\hat f}(v)\,{\hat \psi}(v)\,{\mathrm e}^{2 \pi {\mathrm i}\,t v}\,{\mathrm d}v \\[1ex]
	&= \sum_{k\in \mathbb Z} \frac{1}{L}\,f\big(\tfrac{k}{L}\big)\,\int_{-L/2}^{L/2} {\hat \psi}(v)\,{\mathrm e}^{2 \pi {\mathrm i}\,(t - k/L)\, v}\,{\mathrm d}v %\\[1ex]
	= \sum_{k\in \mathbb Z} f\big(\tfrac{k}{L}\big)\, \frac{1}{L}\,\psi\big(t - \tfrac{k}{L}\big) \,,
	\quad t\in\R \,,
\end{align*}
\new{where summation and integration may be interchanged by the theorem of Fubini--Tonelli, since we have~\eqref{eq:sampledecay} and \mbox{$|\hat\psi(v)|\leq 1$} by definition~\eqref{eq:hatpsi(v)}.
Additionally, note that from~\eqref{eq:hatpsi(v)} and~\eqref{eq:psi} it follows that
$$
|\psi(t)| \leq 2\,\int_0^{L/2} \big| \cos(2\pi v t)\big|\,{\mathrm d}v < L \,, \quad t \in \mathbb R \,.
$$
Hence, we have \mbox{$\frac{1}{L}\,\big| \psi\big(t - \frac{k}{L}\big)\big| < 1$} for all \mbox{$t \in \mathbb R$} and \mbox{$k \in \mathbb Z$}, and consequently the series~\eqref{eq:f=sum} converges absolutely and uniformly on $\mathbb R$ by~\eqref{eq:sampledecay} and the Weierstrass M-test.}
%This completes the proof. 
\hfill\qedsymbol
\medskip

Note that~\eqref{eq:f=sum} is not an interpolating approximation, since in general we have
\begin{align*}
	\frac{1}{L}\,\,\psi\big(t - \tfrac{k}{L}\big)\,\bigg|_{t=\frac nL} \not= \delta_{n,k} \,,\quad
	n,k\in \mathbb Z \,.
\end{align*}
Moreover, since the frequency window function~\mbox{$\hat\psi$} in~\eqref{eq:hatpsi(v)} is compactly supported, the uncertainty principle (cf.~\cite[p.~103, Lemma~2.39]{PPST23}) yields \mbox{$\mathrm{supp}\,\psi = \mathbb R$},
such that~\eqref{eq:f=sum} does not imply localized sampling for any choice of \mbox{$\hat\psi$}.
In other words, the representation~\eqref{eq:f=sum} still requires infinitely many samples \mbox{$f\big(\tfrac{k}{L}\big)$}.
Thus, for practical realizations we need to consider a truncated version of~\eqref{eq:f=sum} and hence for \mbox{$T\in \mathbb N$} we introduce the \mbox{$T$-th} partial sum
\begin{equation}
\label{eq:PTf}
	(P_{\psi, T} f)(t) \coloneqq \sum_{k=-T}^T  f\big(\tfrac{k}{L}\big)\,\frac{1}{L}\,\psi\big(t - \tfrac{k}{L}\big)\,, \quad t\in [-1,\,1] \,.
\end{equation}

Then for the linear frequency window function~\eqref{eq:hatpsilin(v)} we show the following convergence result.

\begin{Theorem}
\label{Thm:convregulpsi}
Let \mbox{$f \in {\mathcal B}_{N/2}(\mathbb R)$} be a bandlimited function with bandwidth \mbox{$\tfrac{N}{2}$}, \mbox{$N \in \N$}, and let \mbox{$L = N\,(1+\lambda)$} with \mbox{$\lambda > 0$} be given.
Assume that the samples \mbox{$f\big(\tfrac{k}{L}\big)$}, \mbox{$k \in \mathbb Z$}, fulfill the condition~\eqref{eq:sampledecay}.
Using oversampling  and regularization
with the linear frequency window function~\eqref{eq:hatpsilin(v)}, the \mbox{$T$-th} partial sums \mbox{$P_{{\mathrm{lin}},T} f$} converge uniformly to~$f$ on \mbox{$[-1,\,1]$} as \mbox{$T \to \infty$}.
For~\mbox{$T > L$} the following estimate holds
\begin{equation}
\label{eq:max|f-PTf|}
	\max_{t\in [-1,\,1]} \big| f(t) - (P_{\mathrm{lin},T} f)(t)\big| \le \sqrt{\frac{2 L}{3}}\,\frac{\new{2}\,(1+\lambda)}{\pi^2 \lambda}\,(\new{T - L})^{-3/2}\,\| f\|_{L^2(\mathbb R)}\,.
\end{equation}
\end{Theorem}

\emph{Proof}. By~\eqref{eq:f=sum} and~\eqref{eq:PTf} we have
\begin{equation*}
	f(t) - (P_{\mathrm{lin},T} f)(t) = \sum_{|k|> T}  f\big(\tfrac{k}{L}\big)\,\frac{1}{L}\,\psi_{\mathrm{lin}}\big(t - \tfrac{k}{L}\big) \,,
\end{equation*}
such that Cauchy--Schwarz inequality implies
\begin{equation}
\label{eq:CSU}
	\big|f(t) - (P_{\mathrm{lin},T} f)(t)\big| \le \bigg( \sum_{|k|> T}  \big|f\big(\tfrac{k}{L}\big)\big|^2\bigg)^{1/2}\,\bigg( \sum_{|k|> T}  \big|\tfrac{1}{L}\,\psi_{\mathrm{lin}}\big(t - \tfrac{k}{L}\big)\big|^2\bigg)^{1/2}\,.
\end{equation}
For \mbox{$f \in {\mathcal B}_{N/2}(\mathbb R)$} the Parseval equality implies
\begin{equation*}
%	\label{eq:Parseval}
	\frac 1L\, \sum_{k \in \mathbb Z} \big|f\big(\tfrac{k}{L}\big)\big|^2 = \|f\|_{L^2(\mathbb R)}^2 \,,
\end{equation*}
see~\cite[Formula (3.8)]{KiPoTa22}, and thereby we have
\begin{equation}
\label{eq:est1}
	\bigg(\sum_{|k| > T} \big|f\big(\tfrac{k}{L}\big)\big|^2\bigg)^{1/2} \le \sqrt L\, \|f\|_{L^2(\mathbb R)}\,.
\end{equation}
It can easily be seen that~\eqref{eq:psilin} satisfies the decay condition
\begin{equation*}
	\big|\tfrac{1}{L}\,\psi_{\mathrm{lin}}(x)\big| \le \frac{2}{LN \lambda \,\pi^2}\,x^{-2}\,, \quad x \in \mathbb R \setminus \{0\}\,,
\end{equation*}
and thereby %for \mbox{$|k| > T > L$} and \mbox{$t \in [-1,\,1]$} we have
\begin{equation*}
	\big|\tfrac{1}{L}\,\psi_{\mathrm{lin}}\big(t - \tfrac{k}{L}\big)\big|^2 \le \frac{4\,(1+\lambda)^2}{\lambda^2 \pi^4}\,(L t - k)^{-4}\,.
\end{equation*}
Thus, for \mbox{$T > L$} and \mbox{$t \in [-1,\,1]$} we obtain
\begin{align*}
%\label{eq:est2}
\bigg( \sum_{|k|> T} \big|\tfrac{1}{L}\,\psi_{\mathrm{lin}}\big(t - \tfrac{k}{L}\big)\big|^2\bigg)^{1/2} &\le \frac{2\,(1+\lambda)}{\lambda \pi^2}\,\bigg( \sum_{|k|> T} (L t - k)^{-4}\bigg)^{1/2} \nonumber \\
&\le \frac{2 \sqrt 2 \,(1+\lambda)}{\lambda \pi^2}\,\bigg( \sum_{k = T+1}^{\infty} (k - L)^{-4}\bigg)^{1/2}\,.
\end{align*}
Using the integral test for convergence of series, we conclude
\begin{align*}
	\sum_{k = T+1}^{\infty} (k - L)^{-4} &
	\le \int_{\new{T}}^{\infty} (t - L)^{-4}\,{\mathrm d}t
	= \tfrac{1}{3}\,(T-L)^{-3} \,,
\end{align*}
which yields
\begin{equation}
\label{eq:est3}
	\bigg( \sum_{k = T+1}^{\infty} (k - L)^{-4}\bigg)^{1/2} \le \tfrac{1}{\sqrt 3}\,(T-L)^{-3/2}\,.
\end{equation}
Therefore, \eqref{eq:CSU}, \eqref{eq:est1}, and~\eqref{eq:est3} imply the estimate~\eqref{eq:max|f-PTf|}.
\hfill\qedsymbol
\medskip

\begin{Example}
\label{ex:visualization_frequency_window}
Next, we visualize the error bound of Theorem~\ref{Thm:convregulpsi}, i.e., for a given function \mbox{$f\in {\mathcal B}_{N/2}(\mathbb R)$} with \mbox{$L= N\,(1+\lambda)$}, \mbox{$\lambda > 0$}, we show that the approximation error satisfies~\eqref{eq:max|f-PTf|}.
For this purpose, the error
\begin{equation}
\label{eq:err_psilin}
	\max_{t\in [-1,\,1]} \big| f(t) - (P_{\mathrm{lin},T} f)(t)\big|
\end{equation}
is estimated by evaluating the given function~$f$ as well as its approximation~\mbox{$P_{\mathrm{lin},T} f$}, cf.~\eqref{eq:PTf}, at equidistant points \mbox{$t_s\in [-1,\,1]$}, \mbox{$s=1,\dots,S$}, with \mbox{$S=10^5$}.
Here we study the function \mbox{$f(t) = \sqrt{N} \,\mathrm{sinc}(N \pi t)$}, \mbox{$t \in \mathbb R$}, such that \mbox{$\|f\|_{L^2(\mathbb R)}=1$}.
We fix \mbox{$N=128$} and consider the error behavior for increasing \mbox{$T \in \mathbb N$}.
More specifically, in this experiment we choose several oversampling parameters \mbox{$\lambda\in\{0.5,1,2\}$} and truncation parameters \new{\mbox{$T = 2^c$} with~\mbox{$c \in \{0,\dots,15\}$}}.
The corresponding results are depicted in Figure~\ref{fig:err_const_frequency}.
Note that the error bound in~\eqref{eq:max|f-PTf|} is only valid for \mbox{$T>L$}.
Therefore, we have additionally marked the point \mbox{$T=L$} for each~$\lambda$ as a vertical dash-dotted line.
It can easily be seen that also the error results are much better when \mbox{$T>L$}.
Note, however, that increasing the oversampling parameter~$\lambda$ requires a much larger truncation parameter~$T$ to obtain errors of the same size.
\begin{figure}[ht]
	\centering
	\captionsetup[subfigure]{justification=centering}
	\begin{subfigure}[t]{0.4\textwidth}
		\includegraphics[width=\textwidth]{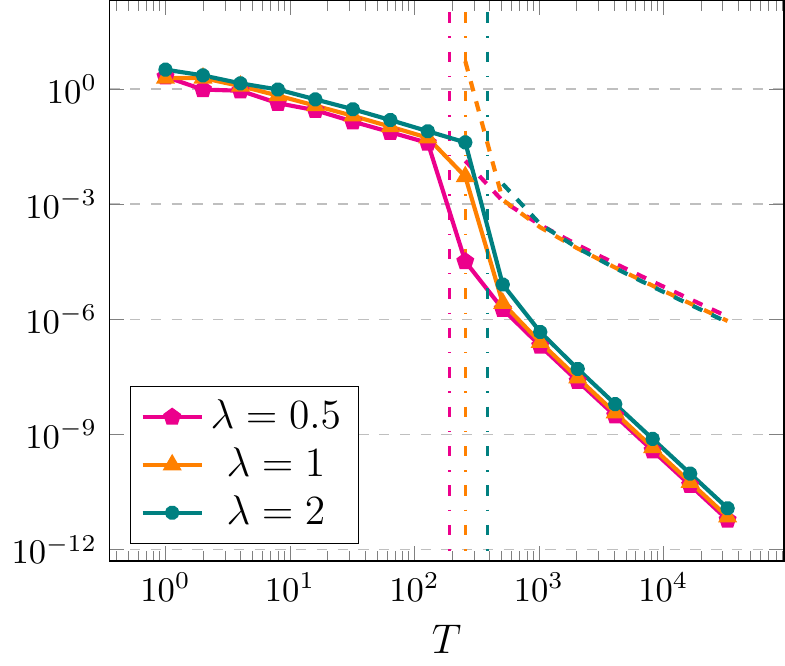}
	\end{subfigure}
	\caption{Maximum approximation error~\eqref{eq:err_psilin} (solid) and error constant~\eqref{eq:max|f-PTf|} (dashed) using the linear frequency window \mbox{$\psi_{\mathrm{lin}}$} from~\eqref{eq:psilin} in~\eqref{eq:PTf} for the function \mbox{$f(t) = \sqrt{N} \,\mathrm{sinc}(N \pi t)$} with \mbox{$N=128$}, \new{\mbox{$T = 2^c$}, \mbox{$c \in \{0,\dots,15\}$}}, and \mbox{$\lambda\in\{0.5,1,2\}$}.
	\label{fig:err_const_frequency}}
\end{figure}
\end{Example}

In order to obtain convergence rates better than the one in Theorem~\ref{Thm:convregulpsi}, one may consider frequency window functions~\eqref{eq:hatpsi(v)} of higher smoothness.

\begin{Example}
\label{Ex:cubicfrequencywindow}
Next, we construct a continuously differentiable frequency window function by polynomial interpolation.
Since the frequency window function~\eqref{eq:hatpsi(v)} is even, it suffices to consider only
\mbox{$\chi \colon \big[\tfrac{N}{2},\,\tfrac{L}{2}\big] \to [0,\,1]$}
at the interval boundaries \mbox{$\tfrac{N}{2}$} and~\mbox{$\tfrac{L}{2}$}.
Clearly, the linear frequency window function~\mbox{$\hat\psi_{\mathrm{lin}}$} in~\eqref{eq:hatpsilin(v)} fulfills
\begin{equation*}
	\lim_{v \to \frac N2} \chi(v) = 1 \,, \quad
	\lim_{v \to \frac L2} \chi(v) = 0 \,.
\end{equation*}
Thus, to obtain a smoother frequency window function, we need to additionally satisfy the first order conditions
\begin{equation*}
	\lim_{v \to \frac N2} \chi^\prime(v) = 0 \,, \quad
	\lim_{v \to \frac L2} \chi^\prime(v) = 0 \,.
\end{equation*}
Then the corresponding interpolation polynomial yields the \emph{cubic frequency window function}
\begin{equation}
\label{eq:hatpsi3(v)}
	{\hat \psi_{\mathrm{cub}}}(v) \coloneqq \left\{ \begin{array}{ll} 1 & \colon |v| \le \tfrac{N}{2}\,, \\
		\frac{16}{\left(L-N\right)^{3}} \, \big(|v| - \frac L2\big)^2 \big(|v|-\frac{3N-L}{4}\big) & \colon \tfrac{N}{2}< |v| < \tfrac{L}{2}\,, \\
		0 & \colon |v| \ge \tfrac{L}{2}\,,
	\end{array} \right.
\end{equation}
see~Figure~\ref{fig:frequency_window_cubic}~(a).
By~\eqref{eq:psi} we see that the inverse Fourier transform of~\eqref{eq:hatpsi3(v)} is given by
\begin{align}
\label{eq:psi3}
	{\psi}_{\mathrm{cub}}(t)
%	&=
%	\frac{12 (\cos(N \pi t) - \cos(L \pi t))}{\pi^4 t^4 \left(L-N\right)^{3}}
%	-\frac{6 (\sin(L \pi t) + \sin(N \pi t))}{\pi^3 t^3 \left(L-N\right)^{2}} \,, \quad
%	t\in\R\setminus \{0\} \\
	&= \new{\frac{L+N}{2}\,\sinc\big(\tfrac{L+N}{2}\,\pi t\big) \cdot \frac{12 \left( \sinc\big(\tfrac{L-N}{2}\,\pi t\big)-\cos\big(\tfrac{L-N}{2}\,\pi t\big) \right)}{\pi^2 t^2 (L-N)^2}} \,, %\quad	t\in\R\setminus \{0\} \,,
\end{align}
for \mbox{$t\in\R\setminus \{0\}$} and~\mbox{$\psi_{\mathrm{cub}}(0) = \frac{L+N}{2}$}, cf. Figure~\ref{fig:frequency_window_cubic}~(b).
\end{Example}

Analogous to Theorem~\ref{Thm:convregulpsi}, the following error estimate can be derived.

\begin{Theorem}
\label{Thm:err_est_psi3}
Let \mbox{$f \in {\mathcal B}_{N/2}(\mathbb R)$} be a bandlimited function with bandwidth \mbox{$\tfrac{N}{2}$}, \mbox{$N \in \N$}, and let \mbox{$L = N\,(1+\lambda)$} with \mbox{$\lambda > 0$} be given.
Assume that the samples \mbox{$f\big(\tfrac{k}{L}\big)$}, \mbox{$k \in \mathbb Z$}, fulfill the condition~\eqref{eq:sampledecay}. Using oversampling and regularization
with the cubic frequency window function~\eqref{eq:hatpsi3(v)}, the \mbox{$T$-th} partial sums \mbox{$P_{{\mathrm{cub}},T} f$} converge uniformly to~$f$ on \mbox{$[-1,\,1]$} as \mbox{$T \to \infty$}.
For \mbox{$T > L$} the following estimate holds
\begin{equation}
\label{eq:err_est_psi3}
	\max_{t\in [-1,\,1]} \big| f(t) - (P_{\mathrm{cub},T} f)(t)\big|
	\le
	\sqrt{\frac{2 L}{5}}\,\frac{24\,(1+\lambda)^2}{\pi^3 \lambda^2}\,(\new{T - L})^{-5/2}\,\| f\|_{L^2(\mathbb R)} \,.
\end{equation}
\end{Theorem}

\begin{Example}
\label{Ex:coswindow}
Another continuously differentiable frequency window function is given in~\cite{Rap96} as the \emph{raised cosine frequency window function}
\begin{equation}
\label{eq:hatpsi_cos(v)}
	{\hat \psi_{\cos}}(v)
	\coloneqq
	\left\{ \begin{array}{ll} 1 & \colon |v| \le \tfrac{N}{2}\,, \\
		\frac 12 + \frac 12 \,\cos \Big(\tfrac{2 |v| - N}{L-N} \pi\Big)  & \colon \tfrac{N}{2}< |v| < \tfrac{L}{2}\,, \\
		0 & \colon |v| \ge \tfrac{L}{2}\,,
	\end{array} \right.
\end{equation}
see~Figure~\ref{fig:frequency_window_cubic}~(a).
By~\eqref{eq:psi} the corresponding function in time domain can be determined as
\begin{align}
\label{eq:psi_cos}
	{\psi}_{\cos}(t)
	&=
	\new{\frac{L+N}{2}\,\sinc\big(\tfrac{L+N}{2}\,\pi t\big) \cdot \frac{\cos\big(\tfrac{L-N}{2}\,\pi t\big)}{1-t^2 (L-N)^2} , \quad
	t\in\R\setminus \big\{\! \pm\tfrac{1}{L-N} \big\}} ,
\end{align}
and \mbox{$\psi_{\cos}\left(\pm\tfrac{1}{L-N}\right) = \tfrac{L-N}{4}\, \cos\big(\tfrac{N\pi}{L-N}\big)$},
see~Figure~\ref{fig:frequency_window_cubic}~(b).
Note that since~\mbox{${\hat \psi_{\cos}}$} in~\eqref{eq:hatpsi_cos(v)} possesses the same regularity as \mbox{${\hat \psi_{\mathrm{cub}}}$} in~\eqref{eq:hatpsi3(v)}, both frequency window functions meet the same error bound~\eqref{eq:err_est_psi3}, cf. Figure~\ref{fig:comp_err}.
\begin{figure}[ht]
	\centering
	\captionsetup[subfigure]{justification=centering}
	\begin{subfigure}[t]{0.4\textwidth}
		\includegraphics[width=\textwidth]{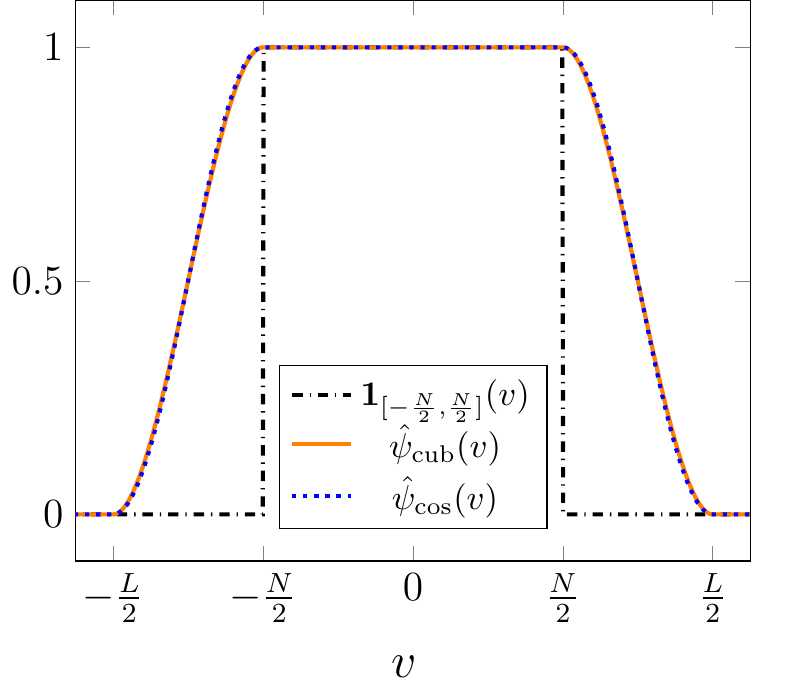}
		\caption{$\hat\psi_{\mathrm{cub}}$ and~$\hat\psi_{\cos}$}
	\end{subfigure}
	\begin{subfigure}[t]{0.4\textwidth}
		\includegraphics[width=\textwidth]{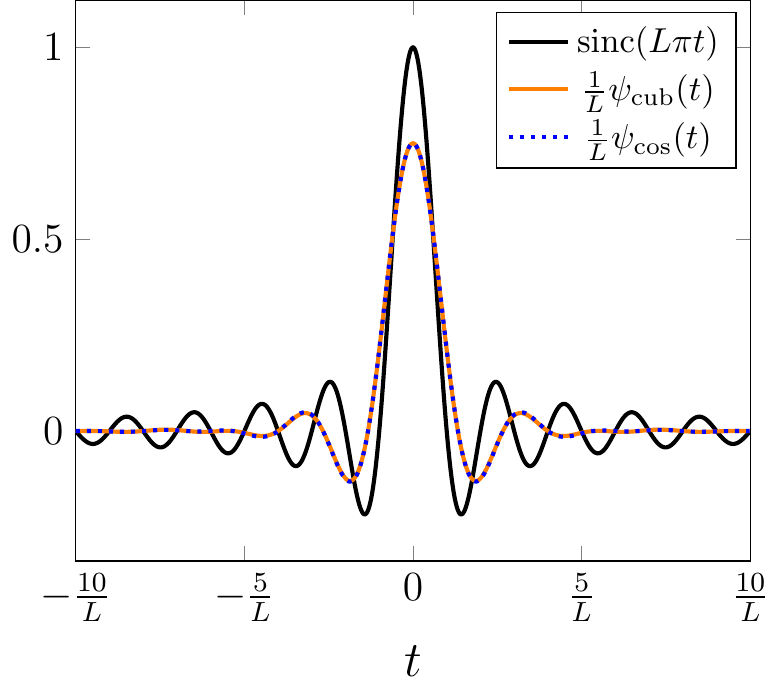}
		\caption{$\frac{1}{L}\,\psi_{\mathrm{cub}}$ and~$\frac{1}{L}\,\psi_{\cos}$}
	\end{subfigure}
	\caption{The frequency window functions~\eqref{eq:hatpsi3(v)} and~\eqref{eq:hatpsi_cos(v)}, and their scaled inverse Fourier transforms.
		\label{fig:frequency_window_cubic}}
\end{figure}
\end{Example}

Note that by~\eqref{eq:psilin} and the convolution property of the Fourier transform, for \mbox{$L > N$} the linear frequency window function~\eqref{eq:hatpsilin(v)} can be written as
\begin{align*}
%\label{eq:hatpsi1(v)_conv}
	{\hat \psi}_{\mathrm{lin}}(v) = \frac{2}{L - N}\,\big({\mathbf 1}_{(N+L)/4} \ast {\mathbf 1}_{({L-N})/4}\big)(v)\,.
\end{align*}
Therefore, instead of determining smooth frequency window functions of the form~\eqref{eq:hatpsi(v)} by means of interpolation as in Example~\ref{Ex:cubicfrequencywindow}, they can also be constructed by convolution, cf.~\cite{Nat86}.

\begin{Lemma}
\label{Lemma:convwithrho}
Let \mbox{$L>N$} be given. Assume that \mbox{$\rho \colon \mathbb R \to [0,\,\infty)$} is an even integrable function with
\mbox{$\mathrm{supp}\,\rho = \big[-\frac{L-N}{4},\,\frac{L-N}{4}\big]$}
and \mbox{$\int_{\mathbb R} \rho(v)\,{\mathrm d}v = 1$}.\\
Then the convolution
\begin{align}
\label{eq:hatpsi_conv}
	{\hat \psi_\mathrm{conv}}(v) = \big({\mathbf 1}_{(N+L)/4} \ast \rho\big)(v)\,, \quad v\in \mathbb R\,,
\end{align}
is a frequency window function of the form~\eqref{eq:hatpsi(v)}.
\end{Lemma}

\emph{Proof.} By assumptions we have
\begin{align}
\label{eq:hatpsi_conv_ext}
	{\hat \psi}_\mathrm{conv}(v)
	= \big({\mathbf 1}_{(N+L)/4} \ast \rho\big)(v)
	= \int_{-(N+L)/4}^{(N+L)/4} \rho(v - w)\,{\mathrm d}w
	= \int_{v-(N+L)/4}^{v+(N+L)/4} \rho(w)\,{\mathrm d}w \,.
\end{align}
Since the convolution of two even functions is again an even function, it suffices to consider~\eqref{eq:hatpsi_conv_ext} only for \mbox{$v \ge 0$}.
For \mbox{$v \in \big[0,\,\frac N2\big]$} we have \mbox{$v - \frac{L+N}{4}\le - \frac{L-N}{4} < \frac{L-N}{4} \le v + \frac{L+N}{4}$} and therefore
\begin{align*}
	{\hat \psi}_\mathrm{conv}(v) =  \int_{-(L-N)/4}^{(L-N)/4} \rho(w)\,{\mathrm d}w = 1 \,.
\end{align*}
For \mbox{$v \in \big[\frac N2,\, \frac L2\big]$} we can write
\begin{align*}
	{\hat \psi}_\mathrm{conv}(v) =  \int_{v-(L+N)/4}^{(L-N)/4} \rho(w)\,{\mathrm d}w \ge 0 \,,
\end{align*}
such that \mbox{${\hat \psi}_\mathrm{conv}\big(\frac N2\big) = 1$}, \mbox{${\hat \psi}_\mathrm{conv}\big(\frac L2\big) = 0$}, and \mbox{${\hat \psi}_\mathrm{conv} \colon \big[\frac N2,\, \frac L2\big] \to [0,\,1]$} is monotonously non-increasing, since \mbox{$\rho(w)\geq 0$} by assumption.
For \mbox{$v \in \big[\frac L2,\,\infty\big)$} we have \mbox{$v - \frac{L+N}{4} \ge \frac{L-N}{4}$}, which implies by assumption \mbox{$\mathrm{supp}\,\rho = \big[-\frac{L-N}{4},\,\frac{L-N}{4}\big]$} that
\begin{align*}
	{\hat \psi}_\mathrm{conv}(v) =  \int_{v-(L+N)/4}^{v+(L+N)/4} \rho(w)\,{\mathrm d}w = 0 \,.
\end{align*}
This completes the proof. \hfill\qedsymbol
\medskip

\new{Given such a frequency window function~\mbox{$\hat\psi_\mathrm{conv}$} as in~\eqref{eq:hatpsi_conv}, its inverse Fourier transform~\eqref{eq:psi} is known by the convolution property of the Fourier transform as
\begin{align}
\label{eq:psiconv}
	\psi_\mathrm{conv}(t) = \frac{N + L}{2} \, \mathrm{sinc}\big(\tfrac{N + L}{2}\,\pi t\big)\, {\check \rho}(t)\,.
\end{align}
Thus, to obtain a suitable window function~\eqref{eq:hatpsi_conv}, we need to assure that the inverse Fourier transform~\mbox{$\check \rho$} of~$\rho$ is explicitly known.
Since~$\rho$ is even by assumption, we have~\mbox{$\check \rho = \hat \rho$} with
\begin{align*}
	\check \rho(t) = \int_{\mathbb R}\rho(v)\,{\mathrm e}^{2 \pi {\mathrm i}\,t v}\,{\mathrm d}v = 2\, \int_0^{\infty} \rho(v)\,\cos(2\pi \,v t)\,{\mathrm d}v\,.
\end{align*}
Note that the convolutional approach~\eqref{eq:hatpsi_conv} has the substantial advantage that the smoothness of~\eqref{eq:hatpsi_conv} is determined by the smoothness of the chosen function~$\rho$.}

\begin{Remark}
\new{The frequency window functions~\mbox{${\hat\psi}_\mathrm{cub}$} in~\eqref{eq:hatpsi3(v)} and~\mbox{${\hat\psi}_{\cos}$} in~\eqref{eq:hatpsi_cos(v)} lack a convolutional representation~\eqref{eq:hatpsi_conv}.
Although the corresponding functions~\eqref{eq:psi3} and~\eqref{eq:psi_cos} in spatial domain are of the form~\eqref{eq:psiconv}, for both frequency windows the Fourier transform of the respective function~$\check\rho$ is only known in the sense of tempered distributions.}
\end{Remark}

\begin{Example}
\label{Ex:conv_Bspline}
For the special choice of
\mbox{$\rho(v) =\frac{2 n}{L - N}\,M_n\big(\frac{2n}{L - N}\,v \big)$} with \mbox{$n \in \mathbb N$}, where $M_n$ is the centered cardinal B-spline of order~$n$, we have
\begin{align*}
	{\check \rho}(t) = \Big(\mathrm{sinc}\big(\tfrac{L-N}{2n}\,\pi t\big)\Big)^n \,.
\end{align*}
Using \mbox{$n=1$} this again yields~\eqref{eq:psilin}, whereas for \mbox{$n=2$} we obtain
\begin{align}
\label{eq:psi_conv2}
	{\psi}_{\mathrm{conv},2}(t) = \frac{N + L}{2} \, \mathrm{sinc}\big(\tfrac{N + L}{2}\,\pi t\big)\, \Big(\mathrm{sinc}\big(\tfrac{L-N}{4}\,\pi t\big)\Big)^2 \,.
\end{align}
Note that the frequency window function~\mbox{$\hat{\psi}_{\mathrm{conv},2}$}, cf.~\eqref{eq:psi_conv2}, possesses the same regularity as~\mbox{$\hat{\psi}_{\mathrm{cub}}$} in~\eqref{eq:hatpsi3(v)} and~\mbox{$\hat{\psi}_{\mathrm{cos}}$} in~\eqref{eq:hatpsi_cos(v)}, and therefore they all meet the same error bound~\eqref{eq:err_est_psi3}, cf. Figure~\ref{fig:comp_err}.
\end{Example}

\begin{Example}
\label{Ex:Natterer}
In \cite{Nat86} the infinitely differentiable function
\begin{align*}
	\rho_{\infty}(v)
	=
	\left\{ \begin{array}{ll} c\,\exp\big(\big[\big(\frac{4v}{L-N}\big)^2 - 1\big]^{-1}\big) \quad &\colon |v| < \frac{L-N}{4}\,,\\
	0 &\colon \text{otherwise}\,, \end{array} \right.
\end{align*}
with the scaling factor
\begin{align*}
	c = \frac{1}{2}\,\bigg( \int_0^{(L-N)/4} \exp\big(\big[\big(\tfrac{4v}{L-N}\big)^2 - 1\big]^{-1}\big)\,{\mathrm d}v\bigg)^{-1} \,.
\end{align*}
is considered.
The corresponding frequency window function~\eqref{eq:hatpsi_conv} is denoted by~\mbox{$\hat\psi_{\infty}$}.
However, since for this function~\mbox{$\rho_{\infty}$} the inverse Fourier transform~\mbox{$\check\rho_{\infty}$} cannot explicitly be stated, the function~\eqref{eq:psiconv} in time domain can only be approximated, which was done by a piecewise rational approximation~\mbox{$\check\rho_{\mathrm{rat}}$} in~\cite{Nat86}.
We remark that because of this additional approximation a numerical decay of the expected rate is doubtful, since the issue of robustness of the corresponding regularized Shannon series remains unclear.
This effect can also be seen in Figure~\ref{fig:comp_err}, where the corresponding frequency window function~\eqref{eq:hatpsi_conv}, denoted by \mbox{$\hat\psi_{\mathrm{rat}}$}, shows similar behavior as the classical Shannon sampling sums~\eqref{eq:ShannonSamplingSum}.

The same comment also applies to \cite{StTa05}, where an infinitely differentiable window function~$\hat\psi$ is aimed for as well.
Since such~$\hat\psi$ with explicit inverse Fourier transform~\eqref{eq:psi} is not known, in~\cite{StTa05} the function~$\psi$  is estimated by some Gabor approximation.
Although an efficient computation scheme via fast Fourier transform (FFT) was introduced in~\cite{StTa06}, the numerical nonrobustness by this approximation seems to be neglected in this work.
%Note in addition that the claimed exponential convergence here does not reflect the maximum error but only some decay on the inside the interval of recovery.
\end{Example}

Finally, we remark that already in \cite[p.~19]{D92} it was stated that a faster decay than for~\mbox{$\hat \psi_{\mathrm{lin}}$} from~\eqref{eq:hatpsilin(v)} can be obtained by choosing~\mbox{$\hat\psi$} in~\eqref{eq:hatpsi(v)} smoother, but at the price of a very large constant.
This can also be seen in Figure~\ref{fig:comp_err}, where the results for the window functions~\mbox{${\hat \psi_{\mathrm{cub}}}$} in~\eqref{eq:hatpsi3(v)}, \mbox{${\hat \psi_{\cos}}$} in~\eqref{eq:hatpsi_cos(v)}, \mbox{${\hat \psi_{\mathrm{conv},2}}$} in~\eqref{eq:psi_conv2}, and \mbox{$\hat\psi_{\mathrm{rat}}$} from Example~\ref{Ex:Natterer} are plotted as well.
For this reason many authors restricted themselves to the linear frequency window function~\mbox{$\hat \psi_{\mathrm{lin}}$} in~\eqref{eq:hatpsilin(v)}.
Furthermore, the numerical results in Figure~\ref{fig:comp_err} encourage the suggestion that in practice only algebraic decay rates are achievable by regularization with a frequency window function.

\section{Regularization with a time window function \label{sec:time_window}}

To preferably obtain better decay rates, we now consider a second regularization technique, namely regularization with a convenient window function in the time domain.
To this end, for~\mbox{$L > N$} and any \mbox{$m \in {\mathbb N} \setminus \{1\}$} with \mbox{$2m \ll L$}, we introduce the set \mbox{$\Phi_{m/L}$}
of all window functions \mbox{$\varphi \colon\,\mathbb R \to [0,\,1]$} with the following properties:\\
\begin{itemize}[leftmargin=*,nosep]
	\item[$\quad\bullet$] Each \mbox{$\varphi \in \Phi_{m/L}$} is supported on \mbox{$\big[-\frac{m}{L},\, \frac{m}{L}\big]$}. Further, $\varphi$ is even and continuous on \mbox{$\big[-\frac{m}{L},\, \frac{m}{L}\big]$}.\\
	\item[$\bullet$] Each \mbox{$\varphi \in \Phi_{m/L}$} restricted on \mbox{$\big[0,\,\frac{m}{L}\big]$} is monotonously non-increasing with \mbox{$\varphi(0) = 1$}.\\
%	%
%	\item[$\bullet$] For each \mbox{$\varphi \in \Phi_{m/L}$}, the Fourier transform
%	%
%	\begin{equation*}
%		{\hat \varphi}(v) = \int_{\mathbb R} \varphi(t)\,{\mathrm e}^{- 2 \pi {\mathrm i}\,v t}\,{\mathrm d}t
%		= 2 \,\int_0^{m/L}  \varphi(t)\,\cos( 2 \pi v t)\,{\mathrm d}t\,, \quad v \in \mathbb R\,,
%	\end{equation*}
%	%
%	is explicitly known.
\end{itemize}
As examples of such window functions we consider the $\sinh$-{\emph{type window function}}
\begin{align}
\label{eq:varphisinh}
	\varphi_{\sinh}(t) \coloneqq
	\begin{cases}
		\frac{1}{\sinh \beta}\, \sinh\Big(\beta\,\sqrt{1-\big(\tfrac{Lt}{m}\big)^2}\,\Big) & \quad \colon t \in \big[-\tfrac{m}{L},\,\tfrac{m}{L}\big] \,, \\
		0 & \quad \colon t \in \mathbb R\setminus \big[-\tfrac{m}{L},\,\tfrac{m}{L}\big]  \,,
	\end{cases}
\end{align}
with \mbox{$\beta = \frac{\pi m\,(L-N)}{L}$},
and the \emph{continuous Kaiser--Bessel window function}
\begin{align}
\label{eq:varphicKB}
	\varphi_{\mathrm{cKB}}(t) \coloneqq
	\begin{cases}
		\frac{1}{I_0(\beta) - 1}\, \bigg( I_0\Big(\beta\,\sqrt{1-\big(\tfrac{Lt}{m}\big)^2}\,\Big) - 1\bigg) & \quad \colon t \in \big[-\tfrac{m}{L},\,\tfrac{m}{L}\big] \,, \\
		0 & \quad \colon t \in \mathbb R\setminus \big[-\tfrac{m}{L},\,\tfrac{m}{L}\big]  \,,
	\end{cases}
\end{align}
with \mbox{$\beta = \frac{\pi m\,(L-N)}{L}$}, where \mbox{$I_0(x)$} denotes the \emph{modified Bessel function of first kind} given by
\begin{equation*}
	I_0(x) \coloneqq \sum_{k=0}^{\infty} \frac{1}{\big((2k)!!\big)^2}\,x^{2k}\,, \quad x\in \mathbb R\,.
\end{equation*}
Both window functions are well-studied in the context of the non\-equispaced fast Fourier transform (NFFT), see e.g. \cite{PT21a} and references therein.

\begin{figure}[ht]
	\centering
	\captionsetup[subfigure]{justification=centering}
	\begin{subfigure}[t]{0.409\textwidth}
		\includegraphics[width=\textwidth]{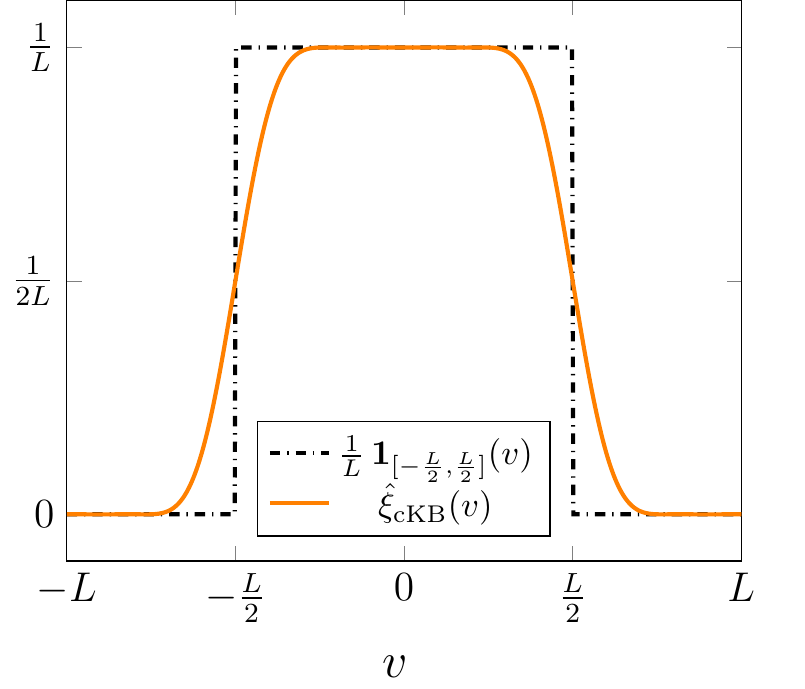}
		\caption{$\new{\hat{\xi}}_{\mathrm{cKB}}$}
	\end{subfigure}
	\begin{subfigure}[t]{0.4\textwidth}
		\includegraphics[width=\textwidth]{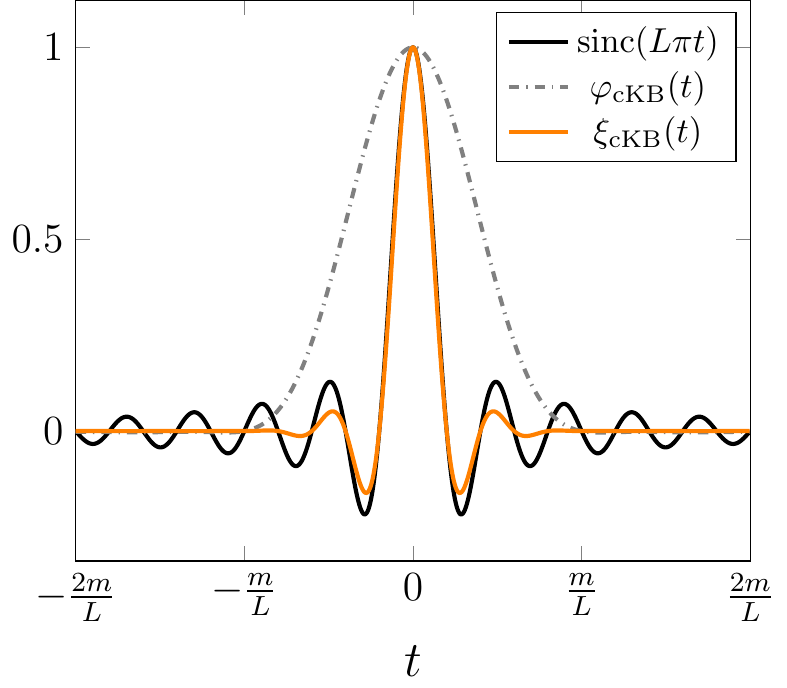}
		\caption{Window function~\eqref{eq:varphicKB}}
	\end{subfigure}
	\caption{The regularized \mbox{$\mathrm{sinc}$} function~\mbox{$\new{\xi}_{\mathrm{cKB}}(t) \coloneqq \mathrm{sinc}(L\pi t) \, \varphi_{\mathrm{cKB}}(t)$} using the continuous Kaiser--Bessel window function~\eqref{eq:varphicKB} and its Fourier transform~\mbox{$\new{\hat{\xi}}_{\mathrm{cKB}}$}.
		\label{fig:cKB_window}}
\end{figure}

A visualization of the continuous Kaiser--Bessel window function~\eqref{eq:varphicKB} as well as the corresponding regularization \mbox{$\new{\xi}_{\mathrm{cKB}}(t) \coloneqq \mathrm{sinc}(L\pi t) \, \varphi_{\mathrm{cKB}}(t)$} of the \mbox{$\mathrm{sinc}$} function can be found in Figure~\ref{fig:cKB_window}.
We remark that in contrast to Figure~\ref{fig:frequency_window} here the function \mbox{$\new{\xi}_{\mathrm{cKB}}$} in time domain is compactly supported and its Fourier transform \mbox{$\new{\hat{\xi}}_{\mathrm{cKB}}$} is
supported on whole~$\mathbb R$, where for the frequency window function~\eqref{eq:hatpsilin(v)} it is vice versa (see \cite[p.~103, Lemma~2.39]{PPST23}).
Note that a visualization for the $\sinh$-type window function~\eqref{eq:varphisinh} would look the same as Figure~\ref{fig:cKB_window}.

Then we \new{approximate} a bandlimited function \mbox{$f \in {\mathcal B}_{N/2}({\mathbb R})$} by the \emph{regularized Shannon sampling formula}
\begin{equation}
\label{eq:regShannon}
	(R_{\varphi,m}f)(t) \coloneqq \sum_{k \in {\mathbb Z}} f\big(\tfrac{k}{L}\big)\,\mathrm{sinc}(L\pi t - \pi k) \, \varphi\big(t - \tfrac{k}{L}\big)\,, \quad t \in {\mathbb R}\,,
\end{equation}
with \mbox{$L \geq N$}.
Since by assumption \mbox{$\mathrm{sinc}(\pi(n-k)) = \delta_{n,k}$} for all \mbox{$n,\,k\in \mathbb Z$} with the Kronecker delta \mbox{$\delta_{n,k}$} and \mbox{$\varphi(0) = 1$},
this procedure is an \emph{interpolating approximation} of~$f$, because
\begin{equation*}
	\mathrm{sinc}(L\pi t - \pi k) \, \varphi\big(t - \tfrac{k}{L}\big)\,\Big|_{t = \frac nL} = \delta_{n,k}\,.
\end{equation*}
Furthermore, the use of the compactly supported window function \mbox{$\varphi \in \Phi_{m/L}$} leads to \emph{localized sampling} of the bandlimited function \mbox{$f \in {\mathcal B}_{N/2}({\mathbb R})$}, i.e.,
the computation of \mbox{$(R_{\varphi,m}f)(t)$} for \mbox{$t \in {\mathbb R} \setminus \tfrac{1}{L}\,{\mathbb Z}$} requires only \mbox{$2m+1$} samples \mbox{$f\big(\frac{k}{L}\big)$}, where \mbox{$k \in
{\mathbb Z}$} fulfills the condition \mbox{$|k - L t| \le m$}. Consequently, for given \mbox{$f \in {\mathcal B}_{N/2}(\mathbb R)$} and \mbox{$L \ge N$}, the reconstruction of~$f$ on the interval \mbox{$[-1,\,1]$} requires
\mbox{$2m+ 2L+1$} samples \mbox{$f\big(\frac{k}{L}\big)$} with \mbox{$k = -m - L,\, \dots,\,m + L$}.
In addition, we again employ oversampling of the bandlimited function \mbox{$f \in {\mathcal B}_{N/2}({\mathbb R})$}, i.e., $f$ is sampled on a finer grid \mbox{$\tfrac{1}{L}\,{\mathbb Z}$} with \mbox{$L >N$}.

This concept of regularized Shannon sampling formulas with localized sampling and oversampling has already been studied by various authors.
A survey of different approaches for window functions can be found in \cite{Q04}, while the prominent Gaussian window function was e.g. considered in \cite{Q03, Q05, SchSt07, TaSuMu07, LZ16}.
Since this Gaussian window function has also been studied in \cite{KiPoTa22}, where superiority of the $\sinh$-type window function~\eqref{eq:varphisinh} was shown, we now focus on time window functions \mbox{$\varphi \in \Phi_{m/L}$}, such as~\eqref{eq:varphisinh} and~\eqref{eq:varphicKB}.

Similar as in \cite{KiPoTa22}, for given \mbox{$f \in {\mathcal B}_{N/2}({\mathbb R})$} and \mbox{$\varphi \in \Phi_{m/L}$} the uniform approximation error
\mbox{$\|f - R_{\varphi,m}f \|_{C_0({\mathbb R})}$}
of the regularized Shannon sampling formula can be estimated as follows.

\begin{Theorem}
\label{Theorem:error-Rmf}
Let \mbox{$f \in {\mathcal B}_{N/2}(\mathbb R)$} be a bandlimited function with bandwidth \mbox{$\tfrac{N}{2}$}, \mbox{$N \in \N$}, and let \mbox{$L > N$}, and \mbox{$m \in {\mathbb N}\setminus \{1\}$} be given.
Further let \mbox{$\varphi \in \Phi_{m/L}$}. \\
Then the regularized Shannon sampling formula~\eqref{eq:regShannon} satisfies
\begin{align*}
%\label{eq:error_const}
	\| f - R_{\varphi,m}f \|_{C_0(\mathbb R)} \le \big( E_1(m,N,L) +  E_2(m,L) \big) \,\|f\|_{L^2(\mathbb R)} \,,
\end{align*}
with the corresponding error constants
\begin{align}
	E_1(m,N,L) &\coloneqq \sqrt N \, \max_{\mathbf v \in [-N/2,\,N/2]} \bigg| 1 - \int_{v - L/2}^{v + L/2} {\hat\varphi}(u)\,{\mathrm d}u \,\bigg|\,, \label{eq:E1} \\
	E_2(m,L) &\coloneqq \frac{\sqrt{2 L}}{\pi\, m}\, \varphi\big(\tfrac{m}{L}\big)\,. \label{eq:E2}
\end{align}
\end{Theorem}

\emph{Proof}. For a proof of Theorem~\ref{Theorem:error-Rmf} see \cite[Theorem~3.2]{KiPoTa22}. \hfill\qedsymbol
\medskip

\new{Note that it is especially beneficial for the estimation of the error constant~\eqref{eq:E1}, if the Fourier transform
\begin{equation}
\label{eq:FTvarphi}
	{\hat \varphi}(v) = \int_{\mathbb R} \varphi(t)\,{\mathrm e}^{- 2 \pi {\mathrm i}\,v t}\,{\mathrm d}t
	= 2 \,\int_0^{m/L}  \varphi(t)\,\cos( 2 \pi v t)\,{\mathrm d}t\,, \quad v \in \mathbb R\,,
\end{equation}
of \mbox{$\varphi \in \Phi_{m/L}$} is explicitly known.
\medskip}

Now we specify the result of Theorem~\ref{Theorem:error-Rmf} for certain window functions.
%First we consider the regularized Shannon sampling formula~\eqref{eq:regShannon} with the $\sinh$-type window function~\eqref{eq:varphisinh}
%which is called
%$\sinh$-\emph{type regularized Shannon sampling formula}
%%
%\begin{equation}
%\label{eq:sinhregShannon}
%	(R_{\sinh,m}f)(t) \coloneqq \sum_{{\mathbf k} \in {\mathbb Z}} f\big(\tfrac{k}{L}\big)\,\mathrm{sinc}(L\pi t - \pi k) \, \varphi_{\sinh}\big(t - \tfrac{k}{L}\big)\,, \quad t \in {\mathbb R}\,.
%\end{equation}
%%
To this end, assume that \mbox{$f \in {\mathcal B}_{N/2}({\mathbb R})$} with \mbox{$N\in \mathbb N$} and \mbox{$L= N\,(1 + \lambda)$}, \mbox{$\lambda > 0$}.
Additionally, we choose \mbox{$m \in \mathbb N \setminus \{1\}$}.
We demonstrate that for the window functions~\eqref{eq:varphisinh} and~\eqref{eq:varphicKB} the approximation error of the regularized Shannon sampling formula~\eqref{eq:regShannon} decreases exponentially with respect to~$m$.
\new{For the $\sinh$-type window function~\eqref{eq:varphisinh} the following result is already known.}

\begin{Theorem}
\label{Thm:error-sinhregShannon}
Let \mbox{$f \in {\mathcal B}_{N/2}(\mathbb R)$} be a bandlimited function with bandwidth \mbox{$\tfrac{N}{2}$}, \mbox{$N \in \N$}, and let \mbox{$L = N\,(1+\lambda)$} with \mbox{$\lambda > 0$} and \mbox{$m \in \mathbb N \setminus \{1\}$} be given.
%Let \mbox{$\varphi_{\sinh}$} be the $\sinh$-type window function~\eqref{eq:varphisinh}.
\\
Then the regularized Shannon sampling formula~\eqref{eq:regShannon} with the $\sinh$-type window function~\eqref{eq:varphisinh} satisfies the error estimate
\begin{equation}
\label{eq:d-error-sinhregShannon}
	\| f - R_{\sinh,m}f \|_{C_0({\mathbb R})} \le \sqrt N\,{\mathrm e}^{-m\pi \lambda/(1+\lambda)}\, \| f \|_{L^2({\mathbb R})}\,.
\end{equation}
\end{Theorem}

\emph{Proof}. For a proof of Theorem~\ref{Thm:error-sinhregShannon} see \cite[Theorem~6.1]{KiPoTa22}. \hfill\qedsymbol
\medskip

Next, we continue with the continuous Kaiser--Bessel window function~\eqref{eq:varphicKB}.

\begin{Theorem}
\label{Thm:error-cKBregShannon}
Let \mbox{$f \in {\mathcal B}_{N/2}(\mathbb R)$} be a bandlimited function with bandwidth \mbox{$\tfrac{N}{2}$}, \mbox{$N \in \N$}, and let \mbox{$L = N\,(1+\lambda)$} with \new{\mbox{$\lambda \geq \frac{1}{m-1}$}} and \mbox{$m \in \mathbb N \setminus \{1\}$} be given.
%Let \mbox{$\varphi_{\mathrm{cKB}}$} be the continuous Kaiser--Bessel window function~\eqref{eq:varphicKB}.
\\
Then the regularized Shannon formula~\eqref{eq:regShannon} with the continuous Kaiser--Bessel window function~\eqref{eq:varphicKB} satisfies the error estimate
\begin{equation}
\label{eq:error-cKBregShannon}
	\| f - R_{\mathrm{cKB},m}f \|_{C_0({\mathbb R})} \le \new{\frac{7\sqrt N\,m\pi\lambda\,(1+\lambda+ 4 m\lambda)}{4\,(1+\lambda)^2} \,{\mathrm e}^{-m\pi\lambda/(1+\lambda)}}\, \| f \|_{L^2({\mathbb R})}\,.
\end{equation}
\end{Theorem}

\begin{figure}[h]
	\centering
	\includegraphics[width=0.4\textwidth]{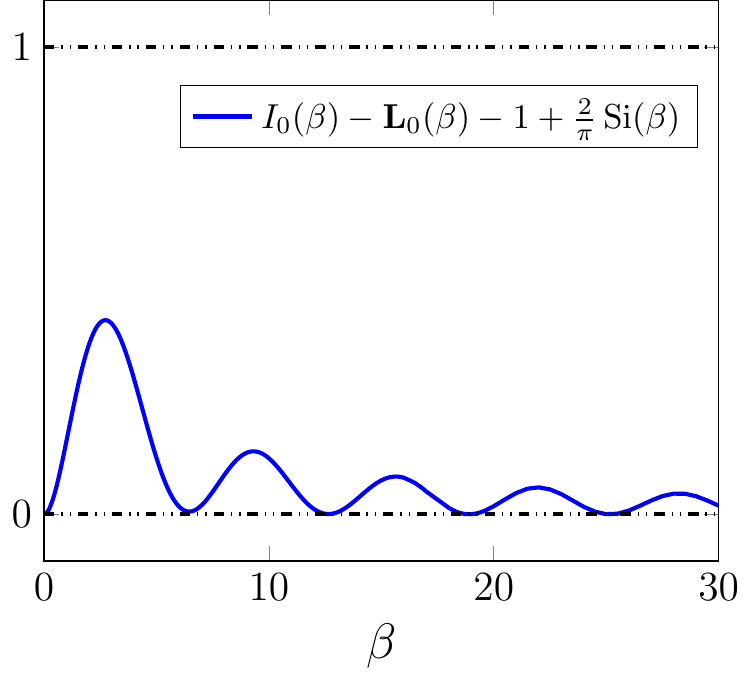}
	\caption{Visualization of \mbox{$\big[ I_0(\beta) - {\textbf L}_0(\beta) - 1 + \frac{2}{\pi}\,\mathrm{Si}(\beta) \big] \in (0,\,1)$} for suitable \mbox{$\beta = \tfrac{m \pi \lambda}{1 + \lambda}$}.
	\label{fig:vis_inequality}}
\end{figure}

\emph{Proof}.
By means of Theorem~\ref{Theorem:error-Rmf} we only have to compute the error constants~\eqref{eq:E1} and~\eqref{eq:E2}.
Note that~\eqref{eq:varphicKB} implies \mbox{$\varphi_{\mathrm{cKB}}\big(\tfrac{m}{L}\big) = 0$}, such that the error constant~\eqref{eq:E2} vanishes.
For computing the error constant~\eqref{eq:E1} we introduce the function \mbox{$\eta \colon\,\big[- \tfrac{N}{2},\,\tfrac{N}{2}\big] \to \mathbb R$} given by
\begin{equation}
\label{eq:eta}
	\eta(v) = 1 - \int_{v-L/2}^{v+L/2} {\hat \varphi}_{\mathrm{cKB}}(u)\,{\mathrm d}u\,.
\end{equation}
As known by \cite[p.~3, 1.1, and p.~95, 18.31]{Ob90}, the Fourier transform of~\eqref{eq:varphicKB} has the form
\begin{align}
\label{eq:FTvarphicKB}
	{\hat \varphi}_{\mathrm{cKB}}(v) = \frac{2 m}{\big(I_0(\beta) - 1\big)\,L}\cdot
	\begin{cases}
		\Big(\frac{\sinh\big(\beta \sqrt{1 - w^2}\,\big)}{\beta \sqrt{1 - w^2}} - \mathrm{sinc}(\beta w)\Big) & \quad \colon |w| < 1 \,, \\
		\big(\mathrm{sinc}\big(\beta \sqrt{w^2-1}\,\big) - \mathrm{sinc}(\beta w)\big) & \quad \colon |w| \ge 1\,,
	\end{cases}
\end{align}
with the scaled frequency \mbox{$w = \frac{2 \pi m}{\beta L}\,v$}.
Thus, substituting \mbox{$w = \frac{2 \pi m}{\beta L}\,u$} in~\eqref{eq:eta} yields
\begin{equation*}
	\eta(v) = 1 - \frac{\beta L}{2 \pi m}\,\int_{-a(-v)}^{a(v)} {\hat \varphi}_{\mathrm{cKB}}\big(\tfrac{\beta L}{2 m \pi}\,w\big)\,{\mathrm d}w
\end{equation*}
with the increasing linear function \mbox{$a(v) \coloneqq \tfrac{2 m \pi}{\beta L}\,\big(v + \tfrac{L}{2}\big)$}.
By the choice of the parameter \mbox{$\beta = \tfrac{m \pi \lambda}{1 + \lambda}$} with \new{\mbox{$\lambda \geq \frac{1}{m-1}$}} we have
\mbox{$a\big(- \tfrac{N}{2}\big) = 1$} and \mbox{$a(v) \ge 1$} for all \mbox{$v \in \big[-\tfrac{N}{2},\,\tfrac{N}{2}\big]$}.
Using~\eqref{eq:FTvarphicKB}, we decompose \mbox{$\eta(v)$} in the form
\begin{equation*}
	\eta(v) = \eta_1(v) - \eta_2(v)\,, \quad v \in \big[-\tfrac{N}{2},\,\tfrac{N}{2}\big] \,,
\end{equation*}
with
\begin{align*}
	\eta_1(v) &= 1 - \frac{\beta}{\pi\,\big(I_0(\beta) - 1\big)}\,\int_{-1}^{1} \bigg(\frac{\sinh\big(\beta \sqrt{1 - w^2}\,\big)}{\beta \sqrt{1 - w^2}} - \mathrm{sinc}(\beta w)\bigg)\,{\mathrm d}w\,, \\[1ex]
	\eta_2(v) &= \frac{\beta}{\pi\,\big(I_0(\beta) - 1\big)}\,\bigg(\int_{-a(-v)}^{-1} + \int_1^{a(v)}\bigg)\Big(\mathrm{sinc}\big(\beta \sqrt{w^2-1}\,\big) - \mathrm{sinc}(\beta w)\Big)\,{\mathrm d}w\,.
\end{align*}
By \cite[3.997--1]{GR80} we have
\begin{align*}
	\int_{-1}^{1} \frac{\sinh\big(\beta \sqrt{1 - w^2}\,\big)}{\beta \sqrt{1 - w^2}}\,{\mathrm d}w &= \frac{2}{\beta}\, \int_0^1 \frac{\sinh\big(\beta \sqrt{1 - w^2}\,\big)}{\sqrt{1 - w^2}}\,{\mathrm d}w \\[1ex]
	&= \frac{2}{\beta}\, \int_0^{\pi/2} \sinh(\beta \cos s)\,{\mathrm d}s = \frac{\pi}{\beta}\,{\textbf L}_0(\beta)\,,
\end{align*}
where \mbox{${\textbf L}_0(x)$} denotes the \emph{modified Struve function} given by (see \cite[12.2.1]{abst})
\begin{equation*}
	{\textbf L}_0(x) \coloneqq \sum_{k=0}^{\infty} \frac{(x/2)^{2k+1}}{\big(\Gamma\big(k + \tfrac{3}{2}\big)\big)^2} = \frac{2 x}{\pi}\,\sum_{k=0}^{\infty}\frac{x^{2k}}{\big((2k+1)!!\big)^2}\,.
\end{equation*}
Note that the function \mbox{$I_0(x) - {\textbf L}_0(x)$} is completely monotonic on \mbox{$[0,\,\infty)$} (see \cite[Theorem~1]{BP14}) and tends to zero as \mbox{$x\to \infty$}.
Applying the \emph{sine integral function}
\begin{equation*}
	\mathrm{Si}(x) \coloneqq \int_0^{x} \frac{\sin w}{w}\,{\mathrm d}w = \int_0^{x} \mathrm{sinc}\, w\,{\mathrm d}w\,, \quad x \in \mathbb R\,,
\end{equation*}
implies
\begin{equation*}
	\int_{-1}^1 \mathrm{sinc}(\beta w)\,{\mathrm d}w = 2\,\int_0^1 \mathrm{sinc}(\beta w)\,{\mathrm d}w = \frac{2}{\beta}\,\mathrm{Si}(\beta)\,.
\end{equation*}
Hence, we obtain
\begin{align*}
	\eta_1(v) &= 1 - \frac{1}{I_0(\beta) - 1}\,\bigg({\textbf L}_0(\beta) - \frac{2}{\pi}\,\mathrm{Si}(\beta)\bigg) %\\
	= \frac{1}{I_0(\beta) - 1}\,\bigg(I_0(\beta) - {\textbf L}_0(\beta) - 1 + \frac{2}{\pi}\,\mathrm{Si}(\beta)\bigg)\,.
\end{align*}
Note that for suitable \mbox{$\beta = \tfrac{m \pi \lambda}{1 + \lambda}$} we find
\mbox{$\big[ I_0(\beta) - {\textbf L}_0(\beta) - 1 + \frac{2}{\pi}\,\mathrm{Si}(\beta) \big] \in (0,\,1)$}, cf. Figure~\ref{fig:vis_inequality}.
In addition, it is known that \mbox{$I_0(x) \geq 1$}, \mbox{$x\in\R$}, such that \mbox{$\eta_1(v) > 0$}.

Now we estimate \mbox{$\eta_2(v)$} for \mbox{$v \in \big[- \tfrac{N}{2},\,\tfrac{N}{2}\big]$} by the triangle inequality as
\begin{equation*}
	|\eta_2(v)| \le \frac{\beta}{\pi\,\big(I_0(\beta) - 1\big)}\,\bigg(\int_{-a(-v)}^{-1} + \int_1^{a(v)}\bigg)\Big|\mathrm{sinc}\big(\beta \sqrt{w^2-1}\,\big) - \mathrm{sinc}(\beta w)\Big|\,{\mathrm d}w\,.
\end{equation*}
By \cite[Lemma~4]{PT21} we have for \mbox{$|w| \ge 1$} that
\begin{equation*}
	\Big|\mathrm{sinc}\big(\beta \sqrt{w^2-1}\,\big) - \mathrm{sinc}(\beta w)\Big| \le \frac{2}{w^2}\,.
\end{equation*}
Thus, we conclude
\begin{equation*}
	|\eta_2(v)| \le \frac{4 \beta}{\pi\,\big(I_0(\beta) - 1\big)}\, \int_1^{\infty} \frac{1}{w^2}\,{\mathrm d}w = \frac{4 \beta}{\pi\,\big(I_0(\beta) - 1\big)}
\end{equation*}
and using Figure~\ref{fig:vis_inequality} we therefore obtain
\begin{align*}
	|\eta(v)| &\le \eta_1(v) + |\eta_2(v)| \le \frac{1}{I_0(\beta)-1}\,\bigg(I_0(\beta) - {\textbf L}_0(\beta) -1 + \frac{2}{\pi}\,\mathrm{Si}(\beta) + \frac{4\beta}{\pi}\bigg)\\
	&\le \frac{1}{I_0(\beta) -1}\,\bigg(1 + \frac{4\beta}{\pi}\bigg)\,.
\end{align*}
\new{By \cite{Ba10} the function~\mbox{$\frac{{\mathrm e}^x}{I_0(x)}$} is strictly decreasing on~\mbox{$[0,\infty)$} and tends to zero as~\mbox{$x \to \infty$}.
Numerical experiments have shown that \mbox{$\frac{{\mathrm e}^{x}}{x\,(I_0(x) - 1)}$} is strictly decreasing on \mbox{$[\pi, \infty)$}, too. By the assumption \mbox{$\lambda \geq \frac{1}{m-1}$} we have
\mbox{$\beta = \frac{\pi m \lambda}{1 + \lambda} \ge \pi$} for all \mbox{$m \in \mathbb N \setminus \{1\}$}. Hence, it follows that
\begin{align*}
	\frac{{\mathrm e}^{\beta}}{\beta\,(I_0(\beta) - 1)} \leq \frac{{\mathrm e}^\pi}{\pi\,(I_0(\pi) - 1)} = 1.644967\ldots < \frac{7}{4}
\end{align*}
and thus we conclude that
\begin{align*}
	\frac{1}{I_0(\beta) - 1}\,\bigg(1 + \frac{4\beta}{\pi}\bigg)
	&< \frac{7\beta}{4}\,\bigg(1 + \frac{4\beta}{\pi}\bigg)\,{\mathrm e}^{-\beta} %\\[1ex]
	= \frac{7\pi m \lambda}{4\,(1+\lambda)}\,\bigg(1 + \frac{4 m \lambda}{1 + \lambda}\bigg)\, {\mathrm e}^{-\pi m \lambda/(1+\lambda)}
\end{align*}
for all  \mbox{$m \in \mathbb N \setminus \{1\}$}.
This completes the proof.}
\hfill\qedsymbol
\medskip

\new{
As seen in Theorem~\ref{Thm:errordata}, if the samples~\mbox{$f\big(\frac{k}{L}\big)$}, \mbox{$k \in \Z$}, of a bandlimited function~\mbox{$f\in {\mathcal B}_{N/2}(\R)$} are not known exactly, i.e., only erroneous samples~\mbox{${\tilde f}_{k} \coloneqq f\big(\frac{k}{L}\big) + \varepsilon_{k}$} with~\mbox{$|\varepsilon_{k}| \leq \varepsilon$}, \mbox{$k\in\Z$}, with~\mbox{$\varepsilon>0$} are known,
the corresponding Shannon sampling series~\eqref{eq:Shannonseries} may differ appreciably from~$f$.
Here we denote the regularized Shannon sampling formula with erroneous samples~\mbox{${\tilde f}_{k}$} by
\begin{align*}
%\label{eq:Rmf_erroneous}
	(R_{\varphi,m}{\tilde f})(t) = \sum_{k\in \Z}  \,\tilde f_k\ \sinc(L \pi t-\pi k) \, \varphi\big(t - \tfrac{k}{L}\big),
	\quad t \in \R.
\end{align*}
Then, in contrast to the Shannon sampling series~\eqref{eq:Shannonseries}, the regularized Shannon sampling formula~\eqref{eq:regShannon} is numerically robust in the worst case analysis,
i.e., the uniform perturbation error \mbox{$\| R_{\varphi,m}{\tilde f} - R_{\varphi,m}f \|_{C_0(\R)}$} is small.
\begin{theorem}
\label{Theorem:robustness-Rmf}
Let \mbox{$f \in {\mathcal B}_{N/2}(\mathbb R)$} be a bandlimited function with bandwidth \mbox{$\tfrac{N}{2}$}, \mbox{$N \in \N$}, and let \mbox{$L > N$}, and \mbox{$m \in {\mathbb N}\setminus \{1\}$} be given.
Further let~\mbox{$\varphi \in \Phi_{m/L}$} as well as~\mbox{${\tilde f}_{k} = f(\frac kL) + \varepsilon_{k}$},
where~\mbox{$|\varepsilon_{k}| \leq \varepsilon$} for all~\mbox{$k\in\Z$}, with~\mbox{$0<\varepsilon\ll 1$}. \\
Then the regularized Shannon sampling sum~\eqref{eq:regShannon} satisfies
\begin{align}
	\| R_{\varphi,m}{\tilde f} - R_{\varphi,m}f \|_{C_0(\R)} &\leq \varepsilon \, \big( 2+L\,\hat\varphi(0) \big) , \label{eq:result_robustness} \\
	\| f - R_{\varphi,m}{\tilde f} \|_{C_0(\R)} &\leq \| f - R_{\varphi,m}{f}\|_{C_0(\R)} + \varepsilon \, \big( 2+L\,\hat\varphi(0) \big) . \notag%\label{eq:robust2}
\end{align}
\end{theorem}
\emph{Proof}. For a proof of Theorem~\ref{Theorem:robustness-Rmf} see \cite[Theorem~3.4]{KiPoTa22}. \hfill\qedsymbol
\medskip}

\new{Note that it is especially beneficial for obtaining explicit error estimates, if the Fourier transform~\eqref{eq:FTvarphi} of \mbox{$\varphi \in \Phi_{m/L}$} is explicitly known.
In the following, we demonstrate that for the window functions~\eqref{eq:varphisinh} and~\eqref{eq:varphicKB} the perturbation error of the regularized Shannon sampling formula~\eqref{eq:regShannon} only grows as~$\mathcal O(\sqrt{m})$.
\begin{theorem}
\label{Theorem:robustness-Rmf-sinh}
Let \mbox{$f \in {\mathcal B}_{N/2}(\mathbb R)$} be a bandlimited function with bandwidth \mbox{$\tfrac{N}{2}$}, \mbox{$N \in \N$}, and let \mbox{$L= N\,(1+\lambda)$} with~\mbox{$\lambda > 0$} and \mbox{$m \in {\mathbb N}\setminus \{1\}$} be given.
Further let~%\mbox{$R_{\mathrm{cKB},m}{\tilde f}$} be as in~\eqref{eq:Rmf_erroneous} with the noisy samples~
\mbox{${\tilde f}_{k} = f\big(\frac{k}{L}\big) + \varepsilon_{k}$},
with~\mbox{$|\varepsilon_{k}| \leq \varepsilon$} for all~\mbox{$k\in\Z$} and~\mbox{$0 <\varepsilon\ll 1$}.\\
Then the regularized Shannon sampling formula~\eqref{eq:regShannon} with the \mbox{$\sinh$-type} window function~\eqref{eq:varphisinh} satisfies the error estimate
\begin{align*}
%\label{eq:result_robustness_sinh}
	\| R_{\sinh,m}{\tilde f} - R_{\sinh,m}f \|_{C_0(\R)}
	\leq
	\varepsilon \left( 2+\sqrt{\frac{2+2\lambda}{\lambda}} \,\frac{1}{1 - \e^{-2 \beta}}\,\sqrt{m} \right).
\end{align*}
\end{theorem}
\emph{Proof}. For a proof of Theorem~\ref{Theorem:robustness-Rmf-sinh} see \cite[Theorem~6.3]{KiPoTa22}. \hfill\qedsymbol
\medskip
\begin{theorem}
\label{Theorem:robustness-Rmf-cKB}
Let \mbox{$f \in {\mathcal B}_{N/2}(\mathbb R)$} be a bandlimited function with bandwidth \mbox{$\tfrac{N}{2}$}, \mbox{$N \in \N$}, and let \mbox{$L= N\,(1+\lambda)$} with~\mbox{$\lambda > 0$} and \mbox{$m \in {\mathbb N}\setminus \{1\}$} be given.
Further let~%\mbox{$R_{\mathrm{cKB},m}{\tilde f}$} be as in~\eqref{eq:Rmf_erroneous} with the noisy samples~
\mbox{${\tilde f}_{k} = f\big(\frac{k}{L}\big) + \varepsilon_{k}$},
with~\mbox{$|\varepsilon_{k}| \leq \varepsilon$} for all~\mbox{$k\in\Z$} and~\mbox{$0 <\varepsilon\ll 1$}.\\
Then the regularized Shannon sampling formula~\eqref{eq:regShannon} with the continuous Kaiser--Bessel window function~\eqref{eq:varphicKB} satisfies the error estimate
\begin{align}
\label{eq:result_robustness_cKB}
	\| R_{\mathrm{cKB},m}{\tilde f} - R_{\mathrm{cKB},m}f \|_{C_0(\R)}
	\leq
	\varepsilon \left( 2 + \sqrt{\frac{2+2\lambda}{\lambda}}\, \sqrt m \right).
\end{align}
\end{theorem}
\emph{Proof}.
By Theorem~\ref{Theorem:robustness-Rmf} we have to compute~\mbox{$\hat{\varphi}_{\mathrm{cKB}}(0)$} for the continuous Kaiser--Bessel window function~\eqref{eq:varphicKB}, which is given by~\eqref{eq:FTvarphicKB} as
\begin{align*}
	\hat{\varphi}_{\mathrm{cKB}}(0)
	&=
	\frac{2 m}{\big(I_0(\beta) - 1\big)\,L}
	\left(\frac{\sinh(\beta)}{\beta} - 1\right)
	=
	\frac{2 m}{L\sqrt{\beta}}
	\cdot \frac{\sinh(\beta)-\beta}{\sqrt{\beta}\,\big(I_0(\beta) - 1\big)} .
\end{align*}
By \cite[9.7.1]{abst} we have
\begin{align*}
	\lim_{x\to\infty} \sqrt{2\pi x} \, \e^{-x} I_0(x) = 1
\end{align*}
and hence
\begin{align*}
\lim_{x\to\infty} \sqrt{2\pi x} \, \e^{-x} \,\big(I_0(x)-1\big) = 1 .
\end{align*}
Moreover, for \mbox{$\beta > 0$} the term \mbox{$\frac{\sinh(\beta)-\beta}{\sqrt{\beta}\,\big(I_0(\beta) - 1\big)}$}
is monotonously increasing with
\begin{align*}
	\lim_{\beta \to \infty} \,\frac{\sinh(\beta)-\beta}{\sqrt{\beta}\,\big(I_0(\beta) - 1\big)}
	=
	\lim_{\beta \to \infty}\, \frac{1-\e^{-2\beta}-2\beta\,\e^{-\beta}}{2\sqrt{\beta}\,\e^{-\beta}\,\big(I_0(\beta) - 1\big)}
	=
	\sqrt{\frac{\pi}{2}}\, ,
\end{align*}
such that~\eqref{eq:result_robustness} and \mbox{$\beta = \frac{\pi m \lambda}{1+\lambda}$} yields the assertion~\eqref{eq:result_robustness_cKB}.
\hfill\qedsymbol
\medskip}

\section{Comparison of the two regularization methods}
\label{sec:comp_err}

Finally, we compare the behavior of the regularization methods presented in Sections~\ref{sec:frequency_window} and~\ref{sec:time_window} to the classical Shannon sampling sums~\eqref{eq:ShannonSamplingSum}.
For a given function \mbox{$f\in {\mathcal B}_{N/2}(\mathbb R)$}
with \mbox{$L= N\,(1+\lambda)$}, \mbox{$\lambda > 0$}, we consider the approxi\-mation errors
\begin{equation}
\label{eq:err_psi}
	\max_{t\in [-1,\,1]} \big| f(t) - (S_{T} f)(t)\big|
	\quad\text{ and }\quad
	\max_{t\in [-1,\,1]} \big| f(t) - (P_{\psi,T} f)(t)\big|
\end{equation}
for \mbox{$\psi \in \{\psi_{\mathrm{lin}},\, \psi_{\mathrm{cub}},\, \psi_{\cos}, \, {\psi}_{\mathrm{conv},2},\, \psi_{\mathrm{rat}}\}$}, cf.~\eqref{eq:psilin}, \eqref{eq:psi3}, \eqref{eq:psi_cos}, \eqref{eq:psi_conv2}, and Example~\ref{Ex:Natterer}, as well as the corresponding error constants~\eqref{eq:max|f-PTf|} and~\eqref{eq:err_est_psi3}.
In addition, we study the approximation error
\begin{equation}
\label{eq:err_phi}
	\max_{t\in [-1,\,1]}| f(t) - (R_{\varphi,m} f)(t) |
\end{equation}
with \mbox{$\varphi \in \{\varphi_{\sinh}, \,\varphi_{\mathrm{cKB}}\}$}, cf.~\eqref{eq:varphisinh} and~\eqref{eq:varphicKB}, and the corresponding error constants~\eqref{eq:d-error-sinhregShannon} and~\eqref{eq:error-cKBregShannon}.
By the definition of the regularized Shannon sampling formula in~\eqref{eq:regShannon} we have
\begin{equation}
\label{eq:Rm}
	(R_{\varphi,m} f)(t) = \sum_{k = -L-m}^{L+m} f(\tfrac k L) \,\new{\xi}(t - \tfrac k L) \,, \quad t \in [-1, 1] \,,
\end{equation}
with the regularized \mbox{$\mathrm{sinc}$} function
\begin{equation}
\label{eq:xi}
	\new{\xi}(t) \coloneqq \mathrm{sinc}(L\pi t) \, \varphi(t) \,.
\end{equation}
Thus, to compare~\eqref{eq:Rm} to \mbox{$S_{T} f$} in~\eqref{eq:ShannonSamplingSum} and~\mbox{$P_{\psi,T} f$} in~\eqref{eq:PTf}, we set \mbox{$T = L+m$}, such that all approximations use the same number of samples.
As in Example~\ref{ex:visualization_frequency_window} the errors~\eqref{eq:err_psi} and~\eqref{eq:err_phi} shall be estimated by evaluating a given function~$f$ and its approximation at equidistant points \mbox{$t_s\in [-1,\,1]$}, \mbox{$s=1,\dots,S,$} with \mbox{$S=10^5$}.
Analogous to \cite[Section~IV,~C]{Ob90}, we choose the function
\begin{equation}
\label{eq:test_func}
	f(t) = \sqrt{\tfrac{4N}{5}} \left[\mathrm{sinc}(N\pi t) + \tfrac 12 \,\mathrm{sinc}(N\pi (t-1))\right]\,, \quad t \in \mathbb R\,,
\end{equation}
with
\mbox{$\|f\|_{L^2(\mathbb R)} = 1$}.
We fix \mbox{$N=256$} and consider several values of \mbox{$m \in \mathbb N \setminus \{1\}$} and \mbox{$\lambda>0$}.

\begin{figure}[t]
	\centering
	\captionsetup[subfigure]{justification=centering}
	\begin{subfigure}[t]{0.62\textwidth}
		\includegraphics[width=\textwidth]{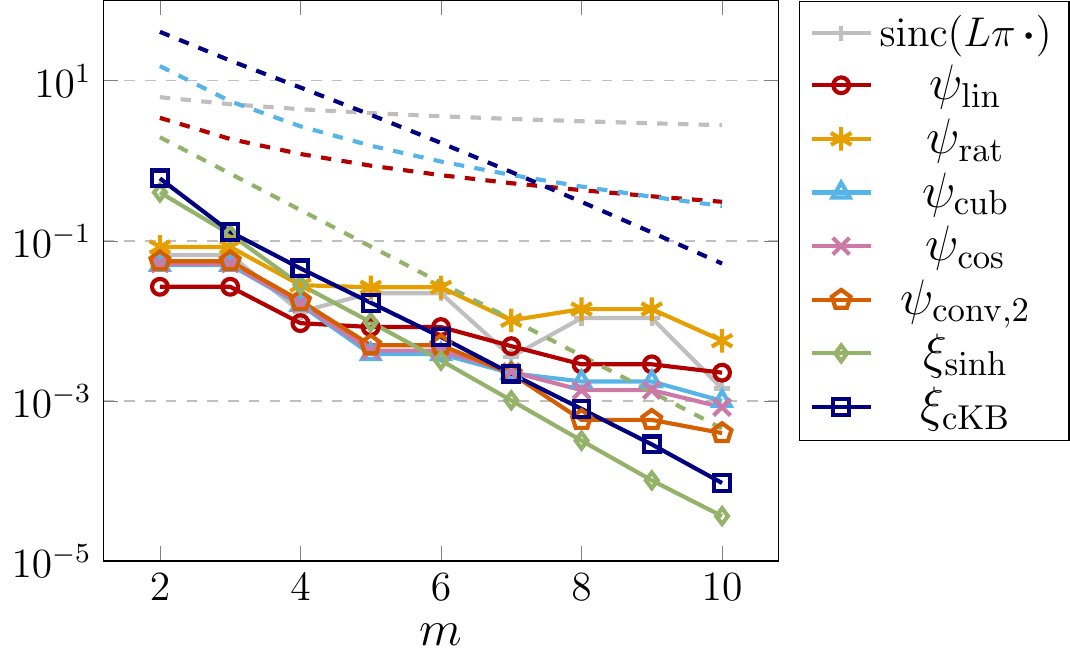}
		\caption{$\lambda=0.5$}
	\end{subfigure}
	\begin{subfigure}[t]{0.62\textwidth}
		\includegraphics[width=\textwidth]{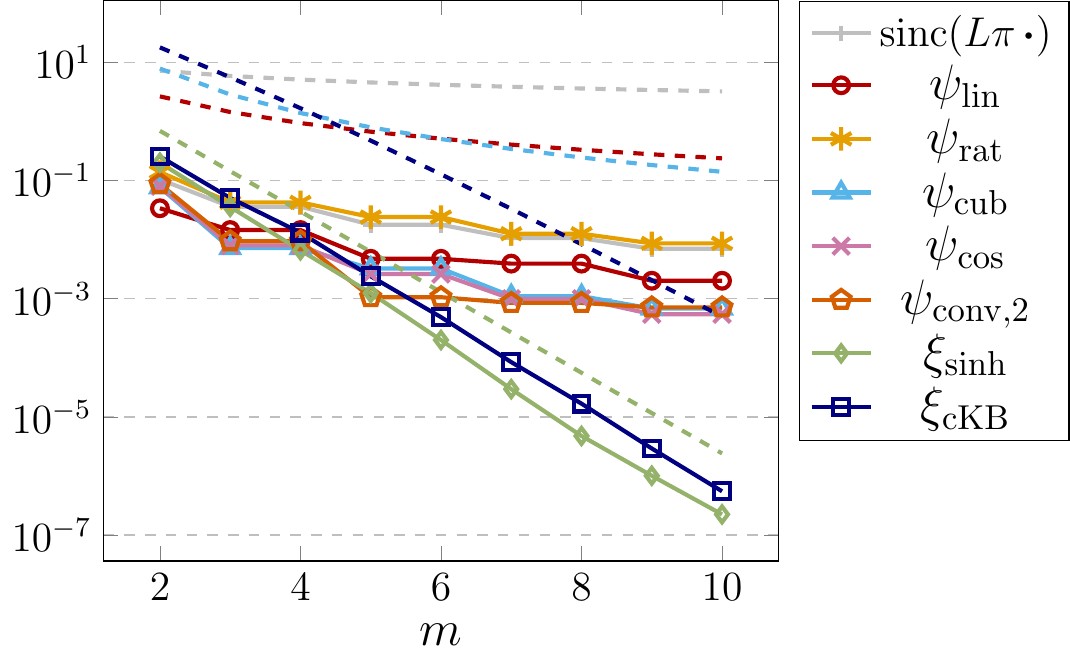}
		\caption{$\lambda=1$}
	\end{subfigure}
	\begin{subfigure}[t]{0.62\textwidth}
		\includegraphics[width=\textwidth]{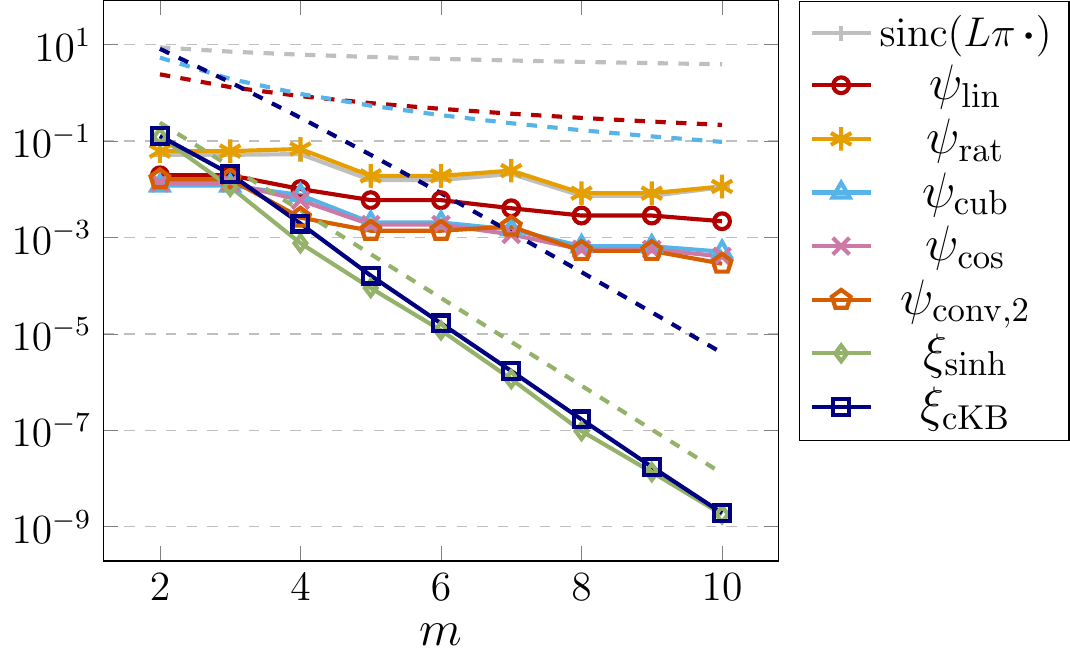}
		\caption{$\lambda=2$}
	\end{subfigure}
	\caption{Maximum approximation error (solid) and error constants (dashed) using \new{classical Shannon sampling sums} compared to regularizations~\eqref{eq:f=sum} with several frequency window functions and regularizations~\eqref{eq:regShannon} with time window functions~\mbox{$\varphi_{\sinh}$} and~\mbox{$\varphi_{\mathrm{cKB}}$}, cf.~\eqref{eq:xi}, for the function~\eqref{eq:test_func} with \mbox{$N=256$}, \mbox{$m\in\{2, 3, \ldots, 10\}$}, and \mbox{$\lambda\in\{0.5,1,2\}$}.
		\label{fig:comp_err}}
\end{figure}

The associated results are displayed in Figure~\ref{fig:comp_err}.
We see that for all window functions the theoretical error behavior perfectly coincides with the numerical outcomes.
\new{In this regard, see also Table~\ref{Table} which summarizes the theoretical results.}
Moreover, it can clearly be seen that for higher oversampling parameter~$\lambda$ and higher truncation parameter~$m$, the error results using~\eqref{eq:regShannon} get much better than the ones using~\eqref{eq:f=sum}, due to the exponential error decay rate shown for~\eqref{eq:regShannon}.

This is to say, our numerical results show that regularization with a time window function performs much better than regularization with a frequency window function, since an exponential decay can (up to now) only be realized using a time window function.
Furthermore, the great importance of an explicit representation of the Fourier transform of the regularizing window function can be seen, cf. Example~\ref{Ex:Natterer}.
\begin{table}
	\centering
	\begin{tabular}[ht]{|c|c|c|}
		\hline\rule{0pt}{3ex}
		window function & error decay rate & see \\[0.5ex]
		\hline\hline\rule{0pt}{3ex}
		$\sinc(L\pi\,\b\cdot)$ & $(T-L)^{-1/2}$ & \cite[Theorem~1]{Ja66} \\
		$\hat\psi_{\mathrm{lin}}$ in~\eqref{eq:hatpsilin(v)} & $(T-L)^{-3/2}$ & Theorem~\ref{Thm:convregulpsi} \\
		$\hat\psi_{\mathrm{cub}}$ in~\eqref{eq:hatpsi3(v)} & $(T-L)^{-5/2}$ & Theorem~\ref{Thm:err_est_psi3} \\
		$\hat\psi_{\mathrm{cos}}$ in~\eqref{eq:hatpsi_cos(v)} & $(T-L)^{-5/2}$ & Example~\ref{Ex:coswindow} \\
		$\hat\psi_{\mathrm{conv},2}$, cf.~\eqref{eq:psi_conv2} & $(T-L)^{-5/2}$ & Example~\ref{Ex:conv_Bspline} \\[0.5ex]
		\hline\rule{0pt}{3ex}
		$\varphi_{\sinh}$ in~\eqref{eq:varphisinh} & ${\mathrm e}^{-m\pi \lambda/(1+\lambda)}$ & Theorem~\ref{Thm:error-sinhregShannon} \\
		$\varphi_{\mathrm{cKB}}$ in~\eqref{eq:varphicKB} & ${\mathrm e}^{-m\pi \lambda/(1+\lambda)}$ & Theorem~\ref{Thm:error-cKBregShannon}
		\\[0.5ex] \hline
	\end{tabular}
	\caption{\new{Summary of the theoretical results on decay rates for the window functions considered in Sections~\ref{sec:frequency_window} and~\ref{sec:time_window}.}
		\label{Table}}
\end{table}

\new{Note that the code files for this and all the other experiments are available on  \url{https://github.com/melaniekircheis/On-numerical-realizations-of-Shannons-sampling-theorem}.}
\medskip

\FloatBarrier

In summary, we found that the regularized Shannon sampling formula with the $\sinh$-type time window function is the best of the considered methods, since this approach is the most accurate, easy to compute, robust in the worst case error, and requires less data (for comparable accuracy) than the classical Shannon sampling sums or the regularization with a frequency window function.

\section*{Acknowledgments}
Melanie Kircheis gratefully acknowledges the support from the BMBF grant 01$\mid$S20053A (project SA$\ell$E).

Moreover, the authors thank the referees and the editor for their very useful suggestions for improvements.

\section*{Conflict of interest statement}
On behalf of all authors, the corresponding author states that there is no conflict of interest.

\end{document}